\documentclass[11pt,a4paper]{article}
\usepackage{color}
\usepackage[colorlinks=true]{hyperref}
\usepackage[margin=1in]{geometry}
\usepackage{fullpage,mathrsfs,graphicx,framed,color,amssymb,amsmath,amsthm}
\usepackage[boxed, noline, ruled, linesnumbered]{algorithm2e}
\usepackage{epsfig, graphicx}
\usepackage{latexsym,amsfonts,amsbsy,amssymb}
\usepackage{amsmath,amsthm}
\usepackage{enumerate}
\usepackage{cleveref}
\usepackage{tikz}
\usepackage{stmaryrd} 
\usepackage{changes} 
\makeatletter 

\@addtoreset{equation}{section}
\makeatother 
\allowdisplaybreaks 
\newtheorem{theorem}{Theorem}[section]
\newtheorem{lemma}{Lemma}[section]

\newtheorem{remark}{Remark}[section]

\allowdisplaybreaks 

\numberwithin{equation}{section}

\begin{document}
\title{Adaptive Neural Network Subspace Method for Solving Partial 
Differential Equations with High Accuracy\footnote{This work was supported by 
the Strategic Priority Research Program of the Chinese Academy of 
Sciences (XDB0640000, XDB0640300),  National Key Laboratory of Computational Physics 
(No. 6142A05230501), National Natural Science Foundations of 
China (NSFC 1233000214), 
National Center for Mathematics and Interdisciplinary Science, CAS.}}
\author{Zhongshuo Lin\footnote{LSEC, Academy of Mathematics and Systems Science,
Chinese Academy of Sciences, No.55, Zhongguancun Donglu, 
Beijing 100190, China, and School of
Mathematical Sciences, University of Chinese Academy
of Sciences, Beijing, 100049 (linzhongshuo@lsec.cc.ac.cn).},\ \ \
Yifan Wang\footnote{School of Mathematical Sciences, Peking University, Beijing 100871, 
China (wangyifan1994@pku.edu.cn).}\ \ \
and \ \ Hehu Xie\footnote{LSEC, Academy of Mathematics and Systems Science,
Chinese Academy of Sciences, No.55, Zhongguancun Donglu, Beijing 100190, 
China, and School of Mathematical Sciences, University of Chinese Academy
of Sciences, Beijing, 100049 (hhxie@lsec.cc.ac.cn).}}
\date{}
\maketitle
\begin{abstract}
Based on neural network and adaptive subspace approximation method, we propose 
a new machine learning method for solving partial differential equations. 
The neural network is adopted to build the basis of the finite dimensional subspace. 
Then the discrete solution is obtained by using the subspace approximation. 
Especially, based on the subspace approximation, a posteriori error estimator 
can be derivated by the hypercircle technique.  This a posteriori error estimator 
can act as the loss function for adaptively refining the parameters of neural network. 


\end{abstract}

\section{Introduction}
Recently, machine learning methods have made a significant impact and 
are becoming increasingly important in science, technology, 
and engineering. So far, there have developed many types of 
machine learning methods for different types of tasks. Among these, the neural network based 
machine learning methods for solving partial differential 
equations (PDEs) attract more and more 
attentions from computational mathematics and scientific 
computing societies. 


It is well known that the machine learning method is built on 
the neural network (NN) and its strong expressing ability. 
The success in the field of image processing motivates 
the large development in many fields such as science, engineering and so on. 
Due to its universal approximation property and strong expressing ability, 
the fully connected neural network (FNN) is the most widely used architecture to 
build the functions for solving PDEs. 
There are several types of well-known FNN-based methods 
such as deep Ritz \cite{EYu}, deep Galerkin method \cite{DGM}, 
PINN \cite{RaissiPerdikarisKarniadakis}, and weak adversarial networks \cite{WAN}
for solving PDEs by designing different loss functions. 
Among these methods, the loss functions always include computing integrations 
for the functions defined by FNN. For example, the loss functions of the deep Ritz method, 
deep Galerkin method and weak adversarial networks method require computing the integrations 
on the computing domain for the functions constructed by FNN. 
Direct numerical integration for the high-dimensional functions always meet 
the ``curse of dimensionality''. Always, the high-dimensional integration  
is computed using the Monte-Carlo method along with some sampling tricks \cite{EYu,HanZhangE}.  
Following this idea, these neural network based machine learning methods always 
use the Monte-Carlo method even on the low-dimensional cases. 
Due to the low convergence rate of the Monte-Carlo method, the solutions obtained 
by the FNN-based numerical methods are challenging to achieve high accuracy and 
stable convergence process. 

Recently, we propose a type of tensor neural network (TNN) and the corresponding
machine learning method to solve high-dimensional problems with high accuracy.
The most important property of TNN is that the corresponding high-dimensional functions
can be easily integrated with high accuracy and high efficiency. Then the deduced machine
learning method can arrive high accuracy for solving high-dimensional problems.
The reason is that the high-dimensional integration of TNN in the loss functions 
can be transformed into one-dimensional integrations which can be computed 
by the classical quadrature schemes with high accuracy.
The TNN based machine learning method has already been used to solve high-dimensional
eigenvalue problems and boundary value problems based on the Ritz type of loss
functions \cite{WangJinXie}. Furthermore,  in \cite{WangXie}, the multi-eigenpairs 
can also been computed with machine learning method which is designed by combining the TNN 
and Rayleigh-Ritz process. Furthermore, with the help of TNN, the 
a posteriori error estimator can be adopted as the loss function of the machine learning method 
for solving high dimensional boundary value problems and eigenvalue 
problems \cite{WangLinLiaoLiuXie}. So far, the TNN based machine learning method has shown 
good ability for solving high dimensional problems.

The idea of finite element methods can give essential motivations for designing 
neural network based machine learning methods for solving PDEs. The aim of this paper is to 
provide the way to exhibit the strong expression (approximation) of NN by solving PDEs 
with high accuracy. 
We will find the machine learning methods with high-accuracy integration possess 
the adaptivity ability as the moving mesh techniques in finite element methods.  
The method in this paper also gives a new way to design machine learning methods 
for solving PDEs and their associated inverse problems. 
The way in this paper shows that the NN based machine learning 
method can achieve high accuracy for solving PDEs. 
The numerical method and analysis are given based on the understanding of 
machine learning from the computational mathematical point of view.  
It is well known that the finite element method is one of the most important and 
popular numerical schemes. The finite element method is built based on 
the finite dimensional subspace approximation for the PDEs. 
Following this idea, the NN based machine learning method for solving PDEs 
can also be designed with the subspace approximation which is built by the NNs. 
Then the process of machine learning method for PDEs is similar to that of the finite 
element method. First, we transform the original PDE to the equivalent variational 
form and then combine with neural network to build the loss function for the machine 
learning method. The quadrature scheme will be considered for computing the 
integrations included in the loss functions. Finally, the discrete loss functions is 
trained (optimized) to deduce the neural network which acts as the approximation 
to the concerned PDE. In this paper, we will consider the second order elliptic 
problems with singularity and interface problems. These problems 
are difficult solved by using the finite element methods. Here, we will investigate 
that the machine learning method can also solve these problems with high accuracy. 

An outline of the paper goes as follows. In Section \ref{Section_Machine_Learning}, 
following the idea of finite element methods, we build the frame for error estimates of 
machine learning method for solving partial differential equations. 
As applications of the understanding in Section \ref{Section_Machine_Learning}, 
we will consider the NN subspace based machine learning method for solving PDEs 
in Section \ref{Section_Subspace}. Section \ref{Section_Interface} is devoted to 
designing the machine learning methods for interface boundary value problems with 
discontinuous coefficients. Some concluding remarks are given in the last section.

\section{Machine learning method for solving PDEs}\label{Section_Machine_Learning} 
In this section, using the NN to replace the finite element space, 
we design a frame of the NN based machine learning for solving PDEs based on the idea of the finite 
element method. Then based on the error estimate theory of finite element method,  
the analysis of the error estimates is also provided for the machine learning method. 

The finite element method for solving PDEs includes three steps. 
The first step is to build the finite element space based on the mesh. Then the second step 
is to build the Ritz and Galerkin form for the concerned PDE on the finite element space. 
The third step is to solve the deduced discrete optimization or linear problem to 
get the finite element solution. The machine learning method for solving PDE can 
also be designed into three steps. The first step is also to build the neural 
network as the trial space. The second step is to build the loss function 
with the Ritz or residual form of NN for the concerned partial differential equations. 
The third step is for the training, which is to optimize the discrete loss function to 
get the neural network solution. From these description, we can find these 
two methods have similar process. 

As we know, the neural network is built by using the compositions of 
linear transforms and nonlinear activation functions. 
From the numerical point of view for solving partial differential equations, 
the neural network provides a new way to build the trial space. 
The successful applications in the field of image processing 
show that the neural network has strong expression ability. 
We also call it strong approximation ability in the numerical point of view. 

In order to show the machine learning method and error estimates 
for the machine learning method, as an example, 
we consider the second order elliptic equation which is defined as: 
Find $u\in H_0^1(\Omega)$ 
such that 
\begin{eqnarray}\label{Elliptic_Problem}
\left\{
\begin{aligned}
-\nabla\cdot\left(\mathcal A\nabla u(x)\right)+\beta u&=f(x),&\ \ \ &x\in\Omega,\\
u(x)&=0,&\ \ \ &x\in\partial\Omega, 
\end{aligned}
\right.
\end{eqnarray}
where $\mathcal A \in\mathbb R^{2\times 2}$ is a symmetric positive definite coefficient matrix and 
the coefficient function $\beta$ is non-negative on the computing domain.

The finite element method is one of the most popular numerical algorithm for solving the second order 
elliptic problem (\ref{Elliptic_Problem}). 
In order to describe the finite element, we define the Ritz form for the 
problem (\ref{Elliptic_Problem}) as follows 
\begin{eqnarray}\label{Ritz_Form}
u = \inf_{v\in H_0^1(\Omega)}\int_\Omega \left(\frac{1}{2}\left(\nabla v\cdot \mathcal A\nabla v + \beta v^2\right)-fv\right)d\Omega.
\end{eqnarray}
From the variational theory, the problem (\ref{Ritz_Form}) is equivalent to the original problem (\ref{Elliptic_Problem}).
The finite element method is defined by restricting the Ritz problem (\ref{Elliptic_Problem}) to the 
finite dimensional subspace which is built by using the piecewise polynomial on the triangulation 
for the computing domain $\Omega$. 
The standard process for finite element method can be described by Algorithm \ref{Process_FEM}.    
\begin{algorithm}
\caption{Process for finite element method}\label{Process_FEM}
\begin{itemize}
\item [1.] Generate the mesh $\mathcal T_h$ for the computing domain 
and build the finite element space $V_h$ 
on the mesh $\mathcal T_h$. The finite element space $V_h$ 
consists of piecewise polynomial basis functions $\varphi_{j,h}$ 
which can be denoted as follows 
\begin{eqnarray}
V_h:=\left\{\varphi_{1,h}, \cdots, \varphi_{N_h,h}\right\},
\end{eqnarray}
where $N_h:={\rm dim} V_h$. 
\item [2.] Build the Ritz form for the problem (\ref{Elliptic_Problem}) 
on the finite element space $V_h$
\begin{eqnarray}\label{Ritz_FEM}
u_h:= \inf_{v_h\in V_h}\int_\Omega\left(\frac{1}{2}
\left(\nabla v_h\cdot\mathcal A \nabla v_h + \beta v_h^2\right)
-fv_h\right)d\Omega,
\end{eqnarray}
where $u_h=\sum_{j=1}^{N_h}c_j \varphi_{j,h}$ denotes the finite element solution. 

\item [3.] Solve the deduced discrete optimization 
problem (\ref{Ritz_FEM}) to obtain the 
coefficient $\{c_j\}_{j=1}^{N_h}$ and the finite element 
solution $u_h=\sum_{j=1}^{N_h}c_j \varphi_{j,h}$.
\end{itemize}
\end{algorithm}
Especially, we should introduce the adaptive mesh redistribution (i.e. Moving Mesh Method) 
in the finite element method \cite{HuangRussell,HuQiaoTang,LiTangZhang}. 
Techniques for the moving mesh method continuously reposition a fixed number 
of mesh nodes and degree of freedom, and so they improve the resolution in 
particular locations of the computational domain.

Similar to the idea of finite element method, the machine learning method for 
problem (\ref{Elliptic_Problem}) can also be defined by restricting 
the Ritz problem (\ref{Ritz_Form}) to the trial space which is built by the neural network. 
Let $V_{\rm NN}$ denote the trial space by neural network and then 
the corresponding restricted Ritz problem is defined as follows 
\begin{eqnarray}\label{Ritz_ML}
u_{\rm NN}^{(1)}:= \inf_{v_N\in V_{\rm NN}}
\int_\Omega\left(\frac{1}{2}\left(\nabla v_N\cdot\mathcal A\nabla v_h 
+\beta v_h^2\right)-fv_N\right)d\Omega,
\end{eqnarray}
where $u_{\rm NN}^{(1)}$ denotes the solution to the problem (\ref{Ritz_ML}). 
While the neural network is not polynomial function, the classical numerical quadratures 
must have integration error. 
In order to do the integration in (\ref{Ritz_ML}), we need to choose the quadrature points $x_i$ and the corresponding 
weights $w_i$, $i=1, \cdots, N$, where $N$ denote the number of quadrature points. 
Based on this quadrature scheme, we can define the computational loss function and the corresponding 
discrete optimization problem as follows
\begin{eqnarray}\label{Ritz_ML_Dist}
u_{\rm NN}^{(2)}:= \inf_{v_N\in V_{\rm NN}}\sum_{i=1}^N\left(\frac{1}{2}
\left(\nabla v_N(x_i)\mathcal A(x_i)\nabla v_N(x_i)
+\beta(x_i)v_N(x_i)^2\right)-f(x_i)v_N(x_i)\right) w_i,
\end{eqnarray}
where $u_{\rm NN}^{(2)}$ denotes the neural network solution to (\ref{Elliptic_Problem}). 

\begin{figure}[htb!]
\centering
\includegraphics[width=13cm, height=3.2cm]{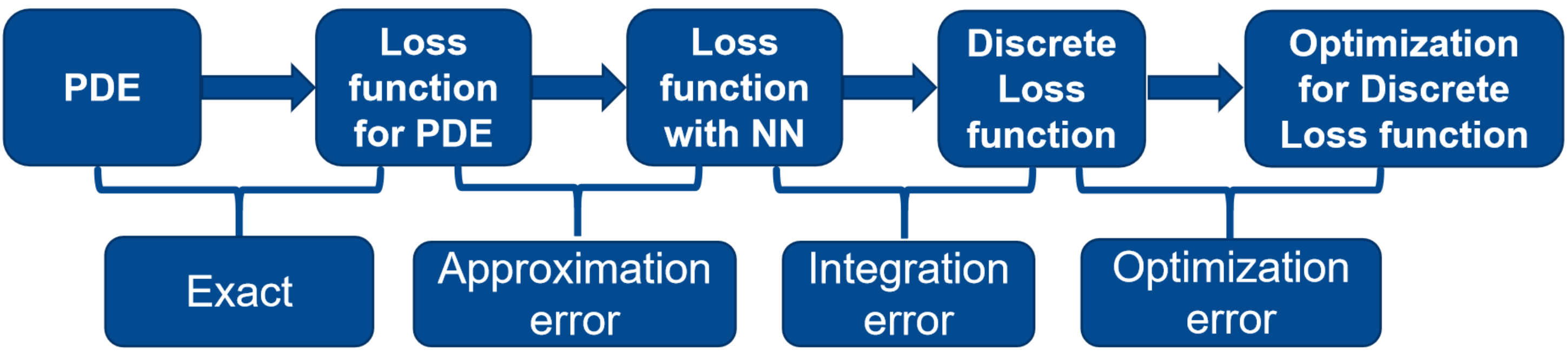}
\caption{Errors of machine learning method for partial differential equations.}\label{Figure_ML_Process}
\end{figure}

The training steps in machine learning is actually doing the optimization for the discrete problem (\ref{Ritz_ML_Dist}). 
We use $u_N^{(3)}$ to denote the approximation solution to $u_N^{(2)}$ after doing sufficient training steps. 
Based on above considerations,  the machine learning method for solving problem (\ref{Elliptic_Problem}) can also 
be described by the process which is given in Algorithm \ref{Process_ML} and Figure \ref{Figure_ML_Process}. 
\begin{algorithm}
\caption{General process for machine learning methd}\label{Process_ML}
\begin{itemize}
\item [1.] Choose some type of architecture to build the neural network $V_{\rm NN}$ 
as the trial space; 

\item [2.] Build the loss function by using the Ritz form and define the corresponding discrete 
problem (\ref{Ritz_ML_Dist}) wiht some type of quadrature scheme. 

\item [3.] Implement sufficient training steps to the  discrete optimize problem (\ref{Ritz_ML_Dist}) 
to obtain the approximation $u_{\rm NN}^{(3)}$ as the final neural network solution. 
\end{itemize}
\end{algorithm}
Let us consider the detailed process of training steps (Step 3) in Algorithm \ref{Process_ML}. 
During the training steps, we adjust the NN parameters to improve the NN solution. 
It should be noted that the number of parameters is fixed across the training process. 
This property is similar to the moving mesh method for finite element scheme. 
Based on this point, we can also understand the NN based machine learning method 
and give the corresponding error analysis following the finite element method. 

For the following analysis, we define the bilinear form $a(\cdot,\cdot)$ as follows
\begin{eqnarray}\label{Bilinear_Form_a}
a(u,v) = \int_\Omega \left(\nabla v \cdot\mathcal A\nabla u + \beta uv\right)d\Omega, \ \ \ 
\forall u\in H_0^1(\Omega),\ \forall v\in H_0^1(\Omega), 
\end{eqnarray}
and then the energy norm as $\|v\|_a=\sqrt{a(v,v)}$.  

It is easy to know that the exact solution $u$ of (\ref{Elliptic_Problem}) 
satisfies the following Galerkin formula 
\begin{eqnarray}\label{Weak_Form}
a(u,v)=(f,v),\ \ \ \forall v\in H_0^1(\Omega),
\end{eqnarray}
where $(\cdot,\cdot)$ denote the $L^2(\Omega)$-inner product 
$$(f,v)=\int_\Omega fvd\Omega.$$

\begin{theorem}\label{Ritz_Approximation_Theorem}
The function $u_{\rm NN}^{(1)}$ defined by (\ref{Ritz_ML}) and the exact solution (\ref{Ritz_Form}) 
satisfy the following error estimate
\begin{eqnarray}\label{Approximation_Error}
\left\|u-u_{\rm NN}^{(1)}\right\|_a = \inf\limits_{v\in V_{\rm NN}}\left\|u-v\right\|_a,
\end{eqnarray}
\end{theorem}
\begin{proof}
From the definitions (\ref{Ritz_Form}) for the exact solution $u$ and 
(\ref{Ritz_ML}) for $u_{\rm NN}^{(1)}$, we have the following optimality
following estimate 
\begin{eqnarray}\label{Equality_1}
\frac{1}{2}\left\|u_{\rm NN}^{(1)}\right\|_a^2-\left(f,u_{\rm NN}^{(1)}\right) 
= \inf\limits_{v\in V_{\rm NN}}\left(\frac{1}{2}\left\|v\right\|_a^2-\left(f, v\right)\right).
\end{eqnarray}
Then based on (\ref{Weak_Form}) and (\ref{Equality_1}), we have the following 
equality for $u$ and $u_{\rm NN}^{(1)}$  
\begin{eqnarray}\label{Equality_2}
\frac{1}{2}\left\|u_{\rm NN}^{(1)}\right\|_a^2-a(u, u_{\rm NN}^{(1)})
+\frac{1}{2}\left\|u\right\|_a^2
= \inf\limits_{v\in V_{\rm NN}}
\left(\frac{1}{2}\left\|v\right\|_a^2-a(u, v)+\frac{1}{2}\left\|u\right\|_a^2\right).
\end{eqnarray}
From (\ref{Equality_2}) and binomial theorem, the following equality hold 
\begin{eqnarray}
\frac{1}{2}\left\|u-u_{\rm NN}^{(1)}\right\|_a^2 
= \inf_{v\in V_{\rm NN}}\frac{1}{2}\left\|u-v\right\|_a^2,
\end{eqnarray}
which leads to the desired result (\ref{Approximation_Error}) and the proof is complete.
\end{proof}

From the definition of the discrete loss function (\ref{Ritz_ML_Dist}) 
and the error estimate for the
numerical integration, we have the following inequality for the error between 
$u_{\rm NN}^{(1)}$ and $u_{\rm NN}^{(2)}$ 
\begin{eqnarray}\label{Integration_Error}
\left\|u_{\rm NN}^{(1)}-u_{\rm NN}^{2}\right\|_a 
\leq \varepsilon_{\rm int}\times {\rm some\ norm\ of}\ 
u_{\rm NN}^{(1)}.
\end{eqnarray}
Then the machine learning method comes to the training steps, 
which is doing the optimization for the 
discrete loss function to obtain the final neural network solution $u_{\rm NN}^{(3)}$. 
For this step,  we can define the optimization error as follows:
\begin{eqnarray}
\left\|u_{\rm  NN}^{(3)}-u_{\rm NN}^{(2)}\right\|_a\leq \varepsilon_{\rm opt}.
\end{eqnarray} 

Based on the analysis, we come to give the final error 
for the neural network solution as follows: 
\begin{eqnarray}\label{Final_Error_Ritz}
\left\|u-u_{\rm NN}^{(3)}\right\|_a
\leq \inf\limits_{v\in V_{\rm NN}}\left\|u-v\right\|_a + \varepsilon_{\rm int} + \varepsilon_{\rm opt}.
\end{eqnarray} 

It is well know that the neural network has strong expression 
ability which leads that the 
approximation error $\inf\limits_{v\in V_{\rm NN}}\left\|u-v\right\|_a$ 
can be small enough with high efficiency. 

From the definitions of loss functions (\ref{Ritz_ML}) and (\ref{Ritz_ML_Dist}), 
there includes integrations which are computed always by using Monte-Carlo 
or sampling methods in normal machine learning methods such as deep-Ritz 
or deep-Galerkin methods. Of course, when $\Omega$ is a low dimensional domain, 
the classical quadrature schemes can also be used to do the integrations with high accuracy. 
A more important fact is that the application of Monte-Carlo or sampling methods always 
leads to the low accuracy of the concerned machine learning methods. For example, 
when using machine learning method with Monte-Carlo integrations for 
solving Laplace eigenvalue problem, the accuracy for the eigenvalue 
approximation is only around 1e-3. 

In the same way, since the strong expression ability of neural network, 
in order to arrive the desired accuracy, the neural network does not need 
to be large. This property leads to that the optimization problem 
(\ref{Ritz_ML_Dist}) is not large neither and not so difficult for solving. 
Based on these consideration, 
the optimization error $\varepsilon_{\rm opt}$ can always be small enough. 

We should also note that the following natural inequality 
\begin{eqnarray}\label{Lower_Bound}
\inf\limits_{v\in V_{\rm NN}}\left\|u-v\right\|_a \leq 
\left\|u-u_{\rm NN}^{(3)}\right\|_a,
\end{eqnarray}
which means that the high accuracy of the neural network  solution for partial 
differential equations leads to the high accuracy of the neural network approximations. 

Then we can draw the conclusion that the integration error always controls 
the machine learning accuracy when the Monte-Carlo or sampling 
method is used for computing the included integrations in loss functions. 
Then, in order to achieve the high accuracy for solving PDEs, it is necessary 
to improve the accuracy of the included integrations in the loss functions. 
For low dimensional cases, it is not so difficult to design numerical quadrature 
schemes with high accuracy. 
Unfortunately, computing high dimensional integration is almost impossible to get high accuracy. 
The Monte Carlo method is the most popular method for computing high dimensional integration. 
But this method always can only obtain low accuracy, which leads to the low accuracy of 
the machine learning method for solving PDEs. Then the essential way to improve the accuracy of 
machine learning method for solving high dimensional PDEs is 
to design the efficient quadrature schemes for high dimensional integration. 
This is the reason to design the tensor neural network and its based machine learning method for 
solving high dimensional PDEs \cite{WangJinXie,WangLinLiaoLiuXie,WangLiaoXie}.  
The important property of TNN is that the high dimensional integration of TNN functions 
can be transformed into one dimensional integrations which can be computed by using the 
classical quadrature schemes with high accuracy. Then the high dimensional integration of TNN 
functions can be computed with high accuracy, which leads to 
the high accuracy of machine learning method for solving high dimensional PDEs.  

The analysis in this section also gives the error estimates of the TNN based machine learning 
method for solving partial differential equations \cite{WangJinXie, WangLinLiaoLiuXie, WangXie}. 
The most important property of TNN is that the high dimensional integration of TNN functions 
can be decomposed into one-dimensional integration which can be computed 
by the classical quadrature schemes with hiah accuracy. This means the integration error 
for TNN based machine learning method is very small which is different from other type  of machine 
learning method for high dimensional problems. This is the essential reason of the high accuracy 
for solving high dimensional problems by using TNN based machine learning method.
 
Based on the similarity to the moving mesh method, the NN based machine learning method should 
have the similar adaptivity ability which is guaranteed by the included training steps. 
Based on these understanding, the machine learning method with high-accuracy integration 
should have high accuracy even for the problems with singularity as that for the moving 
mesh method, which will be investigated in the numerical experiments.

When the PINN or Deep Galerkin scheme is used, the corresponding variational problem takes
the following form 
\begin{eqnarray}\label{Strong_ML}
u = \inf_{v\in H_0^2(\Omega)}\int_{\Omega}\left|-\nabla\cdot\left(\mathcal A\nabla u\right) 
+ \beta u + f\right|^2d\Omega. 
\end{eqnarray}
Based on this variational form, the corresponding restricted problem can be defined as follows 
\begin{eqnarray}\label{Strong_u_N_1}
u_{\rm NN}^{(1)}:=\inf_{v_N\in V_{\rm NN}}\int_\Omega
\left|-\nabla\cdot\left(\mathcal A\nabla v_N \right)+ \beta v_N + f\right|^2d\Omega, 
\end{eqnarray} 
where $u_{\rm NN}^{(1)}$ denotes the solution to the problem (\ref{Strong_ML}).  
In order to do the integration in (\ref{Strong_ML}), we also need to choose the quadrature points $x_i$ 
and the corresponding weights $w_i$, $i=1, \cdots, N$, where $N$ denote the number of quadrature points. 
Based on this quadrature scheme, we can define the computational loss function and the corresponding 
discrete optimization problem as follows
\begin{eqnarray}\label{Strong_ML_Dist}
u_{\rm NN}^{(2)}:=\inf_{v_N\in V_{\rm NN}} \sum_{i=1}^N
\left|-\nabla\cdot\left(\mathcal A(x_i)\nabla v_N(x_i)\right)
+\beta(x_i)v_N(x_i)+f(x_i)\right|^2w_i,
\end{eqnarray}
where $u_{\rm NN}^{(2)}$ denotes a type of neural network solution to (\ref{Elliptic_Problem}). 

Then the training steps in machine learning is actually doing the optimization for the discrete problem (\ref{Strong_ML_Dist}). 
Here, we also use $u_{\rm NN}^{(3)}$ to denote the approximate 
solution to $u_{\rm NN}^{(2)}$ after doing sufficient training steps. 


Based on these consideration, the analysis here can also be extended to 
other neural network based machine learning methods. 

\section{Neural network subspace for solving PDEs}\label{Section_Subspace}
Inspired by the subspace approximation theory of finite element method, we build 
a type of  NN subspace and design the corresponding machine learning method 
for solving low dimensional PDEs. 

As we know, the NN is a type of meshless functions which has strong expression 
ability. The disadvantage for solving PDEs is handling the boundary value condition on the 
complicated domains. This is always the case for low dimensional problems and peoples always 
add penalty terms for the Dirichlet type of boundary value conditions.

In order to describe the subspace based machine learning method, we define the neural network 
subspace as follows
\begin{eqnarray}\label{Subspace_V_p}
V_{p}:={\rm span}\Big\{\varphi_{j}(x;\theta), \ \ j=1, \cdots, p\Big\},
\end{eqnarray}
where $\varphi_{j}(x;\theta)$ denote the $j$-th output of 
the concerned NN $\psi(x,\theta): \mathbb R\rightarrow \mathbb R^p$.  
Based on this subspace $V_p$, we can follow the subspace approximation process to obtain the 
solution $u_p\in V_p$. For this aim, we first define the variational form for 
the problem (\ref{Elliptic_Problem}): Find $u\in V$ such that
\begin{eqnarray}\label{weak_form_hDirichlet}
a(u,v)=(f,v),\ \ \ \forall v\in V, 
\end{eqnarray}
where $V:=H_0^1(\Omega)$. Here and hereafter in this paper, the bilinear form $a(\cdot, \cdot)$ 
is defined as (\ref{Bilinear_Form_a}) and  
the energy norm $\|\cdot\|_a$ is defined as $\|v\|_a=\sqrt{a(v,v)}$.

Based on the $p$-dimensional subspace $V_p$ which satisfies $V_p\subset V$, 
the standard Galerkin approximation scheme for the problem (\ref{weak_form_hDirichlet}) is: 
Find $u_p\in V_p$ such that
\begin{eqnarray}\label{weak_form_hDirichlet_Discrete}
a(u_p,v_p)=(f,v_p),\ \ \ \forall v_p\in V_p.
\end{eqnarray}
Based on the basis of the subspace $V_p$ in (\ref{Subspace_V_p}), we can assemble the 
corresponding stiffness matrix and right hand side term as in the finite element method. 

Then following the idea of the adaptive finite element method or 
moving mesh method, we define the neural network subspace based machine 
learning method for solving problem (\ref{Elliptic_Problem}) as 
in Algorithm \ref{Algorithm_Subspace}. 
\begin{algorithm}[htb!]
\caption{Adaptive neural network subspace method for elliptic problem}\label{Algorithm_Subspace}
\begin{enumerate}
\item Initialization step: Build initial NN $\Psi(x;c^{(0)},\theta^{(0)})$, 
set the loss function $\mathcal L(\Psi(x;c,\theta))$, the maximum training steps $M$, 
learning rate $\gamma>0$, and $\ell=0$.

\item Define the $p$-dimensional space $V_p^{(\ell)}$ as follows
\begin{eqnarray*}
V_{p}^{(\ell)}:={\rm span}\left\{\varphi_{j}(x;\theta^{(\ell)}), \ \ j=1, \cdots, p\right\}.
\end{eqnarray*}
Assemble the stiffness matrix $A^{(\ell)}$ and the right-hand side term 
$B^{(\ell)}$ on $V_{p}^{(\ell)}$ as follows
\begin{eqnarray*}
A_{m,n}^{(\ell)}=a(\varphi_n^{(\ell)},\varphi_m^{(\ell)}),\ \ \ 
B_m^{(\ell)}=(f,\varphi_m^{(\ell)}),\ 1\leq m, n \leq p.
\end{eqnarray*}

\item Solve the following linear equation to obtain the solution 
$c\in\mathbb R^{p\times 1}$ such that 
\begin{eqnarray*}
A^{(\ell)}c=B^{(\ell)}.
\end{eqnarray*}
Update the coefficient parameter as $c^{(\ell+1)}=c$.
Then the Galerkin approximation on the space $V_{p}^{(\ell)}$ for problem 
(\ref{weak_form_hDirichlet}) is $\Psi(x;c^{(\ell+1)},\theta^{(\ell)})$.


\item Update the neural network parameter $\theta^{(\ell)}$ of NN with a 
learning rate $\gamma$
\begin{eqnarray*}
\theta^{(\ell+1)}=\theta^{(\ell)}-\gamma \frac{\partial\mathcal L(\Psi(x;c^{(\ell+1)},\theta^{(\ell)}))}{\partial\theta}.
\end{eqnarray*}

\item Set $\ell=\ell+1$ and go to Step 2 for the next step until $\ell=M$.
\end{enumerate}
\end{algorithm}

\begin{remark}
Sheng and his collaborators developed a subspace-based neural network, 
which includes a subspace layer, to solve partial differential equations (PDEs). 
Their approach involves doing some training steps to the subspace neural 
network to improve the approximation accuracy than the random initial neural network. 
About more information about the random neural networks and their applications to 
solving partial differential equations, please refer to 
\cite{ChenChiEYang, DongLi1,DongLi2,HuangZhuSiew,LiWang,ShangWangSun,SunDongWang}.  
Using the learned subspace, the final solution to the PDE is 
obtained through a least-squares method. In their 
works \cite{DaiFanSheng,LiuXuSheng, XuSheng}, 
they applied this method to solve some types of equations, 
including one-dimensional Helmholtz equations, one-dimensional 
convection equations, one-dimensional parabolic equations, 
two-dimensional Poisson equations, anisotropic diffusion equations, 
and eigenvalue problems, achieving impressive computational accuracy.

Different from the fixed neural network subspace in these methods, the method 
in this paper employs the training steps to improve the accuracy 
by  continuously updating the neural network subspace. 
\end{remark}

\begin{remark}
In \cite{WangLinLiaoLiuXie}, the TNN is adopted to build the NN subspace 
for solving high-dimensional PDEs. Furthermore, based on the finite dimensional 
subspace approximation, the a posteriori error estimators is designed 
for the TNN approximation. Then, following the idea of 
adaptive finite element method, a type of TNN based machine learning method is designed for solving 
high dimensional PDEs with high accuracy, where the a posteriori error estimator act as the loss function 
and the training step is adopted to update the TNN subspace. For more information, please refer to \cite{WangLinLiaoLiuXie}.
\end{remark}

In Algorithm \ref{Algorithm_Subspace}, the NN is updated by minimizing the loss function 
which defines the target of the training steps.  If we use the Ritz type of loss function 
\begin{eqnarray}\label{Ritz_Loss}
\mathcal L(\Psi) = \frac{1}{2}a(\Psi,\Psi)-(f,\Psi),
\end{eqnarray}
as in Theorem \ref{Ritz_Approximation_Theorem}, the training process is to 
minimize the error $\|u-\Psi\|_a$ of the NN approximation $\Psi$ in the energy norm.  
Then the training process is to minimize the error by updating the 
parameters of the fixed NN, which is similar to the moving mesh process of the finite 
element method.

Based on above understanding, the a posteriori error estimators should also be 
adopted to build the loss functions for the machine learning method for solving 
partial differential equations. 
In order to derive the a posteriori error estimate for the subspace approximation $u_p$, 
we introduce the following eigenvalue problem: 
Find $(\lambda,u)\in H_0^1(\Omega)$ such that 
\begin{eqnarray}\label{Eigenvalue_Elliptic}
a(u,v)=\lambda (u,v),\ \ \ \forall v\in V. 
\end{eqnarray}
It is well known that the eigenvalue problem (\ref{Eigenvalue_Elliptic}) has the following 
eigenvalue serires 
\begin{eqnarray*}
0 < \lambda_1 \leq \lambda_2\leq\cdots\leq \lambda_k\leq \cdots\rightarrow+\infty.
\end{eqnarray*}
Then the mimimum eigenvalue $\lambda_{\rm min}=\lambda_1$ and the following inequality holds
\begin{eqnarray}\label{Hilbert_Inequality}
\frac{\|v\|_a^2}{\|v\|_0^2}\geq \lambda_{\rm min}>0,\ \ \ \forall v\in V.
\end{eqnarray}
The following formula of integration by parts acts as the key role to derive the a
posteriori error estimates here. 
\begin{lemma}\label{lemma_Green}
Let $\Omega\subset \mathbb R^d$ be a bounded Lipschitz domain with unit outward normal 
$\mathbf n$ to the boundary $\partial\Omega$.  Then the following Green’s formula holds
\begin{eqnarray}\label{Green_Equality}
\int_\Omega v{\rm div}\mathbf yd\Omega+\int_\Omega \mathbf y\cdot\nabla vd\Omega
=\int_{\partial\Omega}v\mathbf y\cdot\mathbf n ds,\ \ \ \forall v\in H^1(\Omega),
\ \forall\mathbf y\in \mathbf W,
\end{eqnarray}
where $\mathbf W:= \mathbf H({\rm div},\Omega)$.
\end{lemma}
Then the Galerkin approximation $u_p$ satisfies the following a posteriori error estimation.
\begin{theorem}\cite[Theorem 1]{vejchodsky2012complementarity}\label{Theorem_hDirichlet}
Assume $u\in V$ and $u_p\in V_p \subset V$ are the exact solution of 
(\ref{weak_form_hDirichlet}) and subspace approximation $u_p$ 
defined by (\ref{weak_form_hDirichlet_Discrete}), respectively. 
The following upper bound holds
\begin{equation}\label{Upper_Bound_U}
\|u-u_p\|_a \leq \inf_{\mathbf y\in \mathbf W} \eta(u_p,\mathbf y), 
\end{equation}
where $\eta(u_p,\mathbf y)$ is defined as follows
\begin{equation}\label{eta_hDirichlet}
\eta(u_p,\mathbf y):=\left(\|\beta^{-1/2}(f-\beta u_p+{\rm div}(\mathcal A\mathbf y))\|_0^2
+\|\mathcal A^{\frac{1}{2}}(\mathbf y-\nabla u_p)\|_0^2\right)^{\frac{1}{2}}.
\end{equation}
\end{theorem}
This theorem has been proved in the references. For understanding, we also provide the proof
here due to its simplicity.
\begin{proof}
Let us define $v= u-u_p\in V$ in this proof.
Then by inequality (\ref{Hilbert_Inequality}),  Lemma \ref{lemma_Green} and 
Cauchy-Schwarz inequality, the following estimates hold
\begin{eqnarray}\label{Inequality_1_2}
&&a(u-u_p,v)= (f,v)-(\mathcal A\nabla u_p,\nabla v)-(\beta u_p,v)\nonumber\\
&=&(f,v)-(\mathcal A\nabla u_p,\nabla v)-(\beta u_p,v)+(\mathcal A\mathbf y,\nabla v)
+({\rm div}(\mathcal A\mathbf y),v)  - \int_{\partial\Omega}\mathbf n\cdot\mathcal A\mathbf y vds \nonumber\\
&=&(\beta^{-1/2}(f-\beta u_p+{\rm div}(\mathcal A\mathbf y)), \beta^{1/2} v)
+(\mathcal A(\mathbf y-\nabla u_p),\nabla v) \nonumber\\
&\leq& \|\beta^{-1/2}(f-\beta u_p+{\rm div}(\mathcal A\mathbf y))\|_0\|\beta^{1/2}v\|_0
+\|\mathcal A^{\frac{1}{2}}(\mathbf y-\nabla u_p)\|_0
\|\mathcal A^{\frac{1}{2}}\nabla v\|_0\nonumber\\
&\leq& \left(\|\beta^{-1/2}(f-\beta u_p+{\rm div}(\mathcal A\mathbf y))\|_0^2
+\|\mathcal A^{\frac{1}{2}}(\mathbf y-\nabla u_p)\|_0^2\right)^{\frac{1}{2}}
\left(\|\beta^{1/2}v\|_0^2+\|\mathcal A^{\frac{1}{2}}\nabla v\|_0^2\right)^{\frac{1}{2}}\nonumber\\
&\leq& \left(\|\beta^{-1/2}(f-\beta u_p+{\rm div}(\mathcal A\mathbf y))\|_0^2
+\|\mathcal A^{\frac{1}{2}}(\mathbf y-\nabla u_p)\|_0^2\right)^{\frac{1}{2}}
\|v\|_a,\ \ \ \ \forall\mathbf y\in \mathbf W. 
\end{eqnarray}

The desired result (\ref{Upper_Bound_U}) can be deduced easily from (\ref{Inequality_1_2}) and the proof is complete.
\end{proof}

These a posteriori error estimates provide an upper bound of the approximation errors 
in the energy norm. As we know, the numerical methods aim to produce the approximations 
which have the optimal errors in some sense.
This inspires us to use these a posteriori error estimates as the 
loss functions of the NN-based machine learning method
for solving the concerned problems in this paper.

Take the homogeneous Dirichlet type of boundary value problem as an example. 
Denote the Galerkin approximation on space $V_p$ for problem (\ref{weak_form_hDirichlet}) 
as $u_p=u_p(V_p)$.
Now, consider subspaces $V_p\subset V$ with a fixed dimension $p$.
From Theorem \ref{Theorem_hDirichlet}, as long as $V_p\subset V$, the following inequality holds
\begin{eqnarray}\label{ineq_error_eta}
\inf_{\substack{V_p\subset V\\\dim(V_p)=p}}\|u-u_p(V_p)\|_a 
\leq \inf_{\substack{V_p\subset V\\ \dim(V_p)=p}}
\inf_{\mathbf y\in\mathbf W}\eta(u_p(V_p),\mathbf y).
\end{eqnarray}
The error of the best $p$-dimensional Galerkin approximation can be controlled by 
solving the optimization problem on the right-hand side of the inequality (\ref{ineq_error_eta}). 
Based on the discussion in Section \ref{Section_Machine_Learning}, 
we can define the following loss function for the NN-based machine 
learning method at the $\ell$-th step
\begin{eqnarray*}
\mathcal L(u_p)=\inf_{\mathbf y\in\mathbf W}\eta(u_p(V_p^{(\ell)}),\mathbf y).
\end{eqnarray*}
In real implementation, in order to make full use of the parameters of NN and 
improve the computational efficiency, we choose 
$\mathbf y=\nabla u_p \in\mathbf W$. 
Then, we define the loss function as follows
\begin{eqnarray}\label{Loss_Elliptic}
\mathcal L(u_p)=\eta\left(u_p,\nabla u_p\right) 
=\left\|\beta^{-1/2}\left(f-\beta u_p+ \nabla\cdot(\mathcal A\nabla u_p)\right)\right\|_0^2.
\end{eqnarray}

In order to understand the loss function (\ref{Loss_Elliptic}), 
we come to introduce the equality of the errors of the problem \eqref{weak_form_hDirichlet}  
and its dual problem. 
From \cite{vejchodsky2012complementarity}, the minimization problem in (\ref{Upper_Bound_U}) 
is equivalent to another variational problem, which is called the dual problem 
to \eqref{weak_form_hDirichlet} here, that is: Find $\mathbf y^*\in \mathbf W$ such that
\begin{equation}\label{Dual_Problem}
a^*(\mathbf y^*,\mathbf z)= \mathcal F^*(\mathbf z),
\ \ \ \ \forall \mathbf z\in \mathbf W,
\end{equation}
where
\begin{equation*}
a^*(\mathbf y^*,\mathbf z)=\int_{\Omega}\big(\beta^{-1}{\rm div}(\mathcal A\mathbf y^*) 
{\rm div}(\mathcal A\mathbf z)
+\mathbf y^*\cdot\mathcal A\mathbf z\big)d\Omega, \ \  \
\mathcal F^*(\mathbf z)=-\int_{\Omega}\beta^{-1}f{\rm div}(\mathcal A\mathbf z)d\Omega.
\end{equation*}
It is obvious that the bilinear form $a^*(\cdot,\cdot)$ induces an inner product in $\mathbf W$ and the corresponding norm
is $\|\mathbf z\|_{*}=\sqrt{a^*(\mathbf z,\mathbf z)}$. Thanks to the Riesz representation theorem, 
the dual problem \eqref{Dual_Problem} has unique solution.

\begin{theorem}
Let $\beta>0$ and let $u \in V$ be the weak solution of (\ref{weak_form_hDirichlet}). 
Then $\mathbf{y}^*=\boldsymbol{\nabla} u$ solves the dual problem (\ref{Dual_Problem}) 
and $\eta\left(u_h, \mathbf{y}^*\right)=\left\|u-u_h\right\|$.
\end{theorem}

\begin{proof} 
It is easy to know that $\boldsymbol{\nabla} u \in \mathbf{W}$ 
and the exact solution $u$ satisfies the following equality 
$$
\int_{\Omega}\left(-\operatorname{div}(\mathcal A \nabla u)
+\beta u-f\right) v d\Omega=0 \quad \forall v \in V.
$$
Thus $-\operatorname{div} (\mathcal A \boldsymbol{\nabla} u)+\beta u-f=0$ 
a.e. in $\Omega$. Consequently,
$$
\int_{\Omega} \operatorname{div}(\mathcal A\boldsymbol{\nabla} u) \operatorname{div}(\mathcal A\mathbf{w})\mathrm{d}\Omega
-\int_{\Omega} \beta u \operatorname{div}(\mathcal A\mathbf{w})\mathrm{d}\Omega
=-\int_{\Omega} f \operatorname{div}(\mathcal A\mathbf{w})\mathrm{d}\Omega, 
\quad \forall \mathbf{w} \in \mathbf{W}.
$$
Since $$-\int_{\Omega} \beta u \operatorname{div}(\mathcal A\mathbf{w})\mathrm{d}\Omega 
=\int_{\Omega} \beta \boldsymbol{\nabla} u \cdot(\mathcal A\mathbf{w})\mathrm{d}\Omega,
\ \ \ \forall \mathbf{w} \in \mathbf{W},$$ 
we conclude that $\boldsymbol{\nabla} u$ solves (\ref{Dual_Problem}). 
The equality $\eta\left(u_p, \boldsymbol{\nabla} u\right)=\left\|u-u_p\right\|_a$ 
follows immediately from (\ref{eta_hDirichlet}), 
because $f+\operatorname{div}(\mathcal A\boldsymbol{\nabla} u)=$ $\beta u$ a.e. in $\Omega$.
\end{proof}

More importantly, the quantity $\eta(w,\mathbf y)$ has the following complementarity 
result for any $w\in V$ and $\mathbf y\in\mathbf W$.
\begin{lemma}(\cite[Corollary 1]{vejchodsky2012complementarity})\label{Optimization_Property_Lemma}
Assume $\psi$ and $\mathbf y$ are any approximations to the exact solutions $u$ and $\mathbf y^*$
of the problem \eqref{weak_form_hDirichlet} and its dual problem \eqref{Dual_Problem}, respectively.
Then the following equality holds
\begin{equation}\label{Optimization_Property}
\|u-\psi\|_a^2+\|\mathbf y^*-\mathbf y\|_*^2=\eta^2(w,\mathbf y).
\end{equation}
\end{lemma}

\begin{proof}
Putting $\mathbf{w}=\mathbf{y}^*-\mathbf{y}$ and using (\ref{eta_hDirichlet}), 
we may directly compute
\begin{eqnarray}\label{Equality_3}
\eta^2\left(\psi, \mathbf{y}\right) &=&\eta^2\left(\psi, \mathbf{y}^*-\mathbf{w}\right)
=\left\|\beta^{-1/2}\left(f-\beta\psi   
+\operatorname{div}(\mathcal A\mathbf{y}^*)\right)\right\|_0^2 \nonumber\\
&& -2 \int_{\Omega} \beta^{-1}\left(f-\beta \psi
+\operatorname{div}(\mathcal A \mathbf{y}^*)\right) \operatorname{div}(\mathcal A \mathbf{w}) 
\mathrm{d}\Omega 
+\left\|\beta^{-1/2} \operatorname{div}(\mathcal A\mathbf{w})\right\|_0^2 \nonumber\\
&& +\left\|\mathcal A^{1/2}(\mathbf{y}^*-\boldsymbol{\nabla} \psi)\right\|_0^2
-2 \int_{\Omega}\left(\mathbf{y}^*-\boldsymbol{\nabla} \psi\right)
\cdot \mathcal A \mathbf{w} \mathrm{d}\Omega
+\|\mathcal A^{1/2}\mathbf{w}\|_0^2.
\end{eqnarray}
Since $\mathbf{y}^*$ solves (\ref{Dual_Problem}) and due to (\ref{Green_Equality}), 
we have the following orthogonality relation
\begin{eqnarray}\label{Equality_4}
\int_{\Omega} \beta^{-1}\left(f-\beta \psi 
+\operatorname{div}(\mathcal A\mathbf{y}^*)\right) \operatorname{div} (\mathcal A \mathbf{w}) \mathrm{d} \Omega+\int_{\Omega}\left(\mathbf{y}^*-\boldsymbol{\nabla} \psi\right) 
\cdot \mathcal A\mathbf{w} \mathrm{d}\Omega=0,
\quad \forall \mathbf{w} \in \mathbf{W}.
\end{eqnarray}
The proof is finished by substitution of (\ref{Equality_4}) into (\ref{Equality_3}).
\end{proof}
Based on Lemma \ref{Optimization_Property_Lemma}, the loss function in (\ref{Loss_Elliptic})
actually measures the error $\|u-u_p\|_a^2$ and $\|\mathbf y^*-\nabla u_p\|_*^2$, i.e., 
\begin{eqnarray}\label{Loss_Equality}
\mathcal L(u_p)=\|u-u_p\|_a^2 + \|\mathbf y^*-\nabla u_p\|_*^2.
\end{eqnarray}
Furthermore, since $\mathbf y^*=\boldsymbol\nabla u$, then we have 
\begin{eqnarray}\label{Loss_Equality_2}
\mathcal L(u_p)=\|u-u_p\|_a^2 + \|\boldsymbol\nabla u-\boldsymbol\nabla u_p\|_*^2.
\end{eqnarray}

To make the loss function work also for the case $\beta=0$, 
we modify the loss function (\ref{Loss_Elliptic}) to the following form 
\begin{eqnarray}\label{Loss_Elliptic_New}
\mathcal L(u_p)=\eta\left(u_p,\nabla u_p\right) 
=\left\|f-\beta u_p+ \nabla\cdot(\mathcal A\nabla u_p)\right\|_0^2.
\end{eqnarray}

In order to show that the machine learning method has the ability to solve the 
singular problem, we consider solving the Laplace problem on the unit square 
with singularity: Find $u\in H_0^1(\Omega)$ such that 
\begin{eqnarray}\label{Laplace_Singular}
\left\{
\begin{array}{rcl}
-\Delta u&=& f,\ \ \ {\rm in}\ \Omega,\\
u&=&0,\ \ \ {\rm on}\ \partial\Omega,
\end{array}
\right.
\end{eqnarray}
where $\Omega=(0,1)\times (0,1)$ and the right hand side term $f$ 
is chosen  as 
\[
\begin{aligned}
f =& - x_2(x_2 - 1)\left [ \left (x_1 - \frac{1}{2}\right )^{-\frac{2}{3}}
\left (\frac{70}{9}x_1 \left (x_1 - 1 \right ) + \frac{11}{6} \right ) 
+ 2\left (x_2 - \frac{1}{2}\right )^{\frac{4}{3}}\right ]\\
&- x_1(x_1 - 1)\left[\left (x_2 - \frac{1}{2}\right )^{-\frac{2}{3}}
\left(\frac{70}{9}x_2\left(x_2 - 1 \right) 
+ \frac{11}{6}\right) + 2\left(x_1 - \frac{1}{2}\right )^{\frac{4}{3}}\right],
\end{aligned}
\]
such that the exact solution is 
\begin{eqnarray}
u = x_1(1-x_1)x_2(1-x_2)\left(\left(x_1-\frac{1}{2}\right)^{\frac{4}{3}}
+\left(x_2-\frac{1}{2}\right)^{\frac{4}{3}}\right).
\end{eqnarray}
It ie easy to find that there exist ``line'' singularity along 
the lines $x_1=\frac{1}{2}$ and $x_2=\frac{1}{2}$. 
This problem is always used for testing the performance of the 
adaptive finite element methods \cite{mitchell2013collection}. Here, the NN subspace based machine learning method 
is adopted to solve the problem (\ref{Laplace_Singular}). 
For this aim, we build the trial function $u_N$ as follows
\begin{eqnarray}\label{Trial_Function_1}
u_N=x_1(1-x_1)x_2(1-x_2)\Psi(x,\theta).
\end{eqnarray}
Then the trial function $u_N$ is guaranteed to satisfy the 
Dirichlet boundary condition. 

We define the Ritz type of loss function 
for (\ref{Laplace_Singular}) as follows 
\begin{eqnarray}\label{Ritz_Laplace_Singular}
\mathcal L(u_N) = \frac{1}{2}\int_\Omega |\nabla u_N|^2d\Omega 
- \int_\Omega fu_Nd\Omega. 
\end{eqnarray}

In Section~\ref{Section_Machine_Learning}, we provide the error estimate \eqref{Final_Error_Ritz} 
for the neural network solution obtained by the machine learning algorithm. 
From this, it can be observed that the integration error is often the 
main source of error in the machine learning algorithm. 
Therefore, to achieve high-accuracy numerical solutions with machine learning, 
it is essential to ensure the high accuracy of the loss function computation.

In the loss function \eqref{Ritz_Laplace_Singular}, 
the term $\frac{1}{2}\|\nabla u_N\|_0^2$ can be calculated 
using the Gauss-Legendre quadrature. 
However, the singularity of the right-hand side term $f$ 
poses difficulties in achieving high accuracy for the term $(f,u_N)$. 
To handle the singularity of the integral factors 
$(x_1 - \frac{1}{2})^{-2/3}$ and $(x_2 - \frac{1}{2})^{-2/3}$, 
we divide the integration domain into subdomains. 
By applying coordinate transformation, the integral factors are treated 
as weight functions, and the Gauss-Jacobi quadrature is employed 
for computing the integration.

We now briefly introduce the Gauss-Jacobi quadrature which is commonly 
used to compute weighted integrals over $[-1,1]$ 
with the weight function $(1-x)^\alpha(1+x)^\beta$, 
where $\alpha > -1$ and $\beta > -1$. 
Specifically, using $n$ nodes $\{x_i\}$ and weights $\{w_i\}$, 
the Gauss-Jacobi quadrature formula is given by
\[
\int_{-1}^{1}f(x)(1 - x)^\alpha(1+x)^\beta dx \approx \sum_{i=1}^{n}w_if(x_i),
\]
where the $n$ nodes are the roots of the $n$-th Jacobi polynomial 
$P_n^{(\alpha, \beta)}(x)$, and the weights can be computed using the formula:
\[
    w_i = \frac{2^{\alpha+\beta+1}\Gamma(n+\alpha+1)
    \Gamma(n+\beta+1)}{n!
    \Gamma(n+\alpha+\beta+1)(1-x_i^2)[P_n^{(\alpha,\beta)^\prime}(x_i)]^2}.
\]
The Gauss-Jacobi quadrature formula with $n$ nodes is exact 
for polynomials of degree up to $2n-1$. Here, we denote 
the Gauss-Jacobi quadrature formula with the weight function 
$(1-x)^\alpha(1+x)^\beta$ as the $(\alpha,\beta)$-Gauss-Jacobi quadrature. 
For a detailed introduction to the Gauss-Jacobi quadrature, 
please refer to \cite{shen2011spectral}.

Next, we describe the computation for the term $(f,u_N)$.  
Here, we omit the notation $\theta$ in $u_N$ and express it simply as $u_N(x_1, x_2)$:
\[   
\begin{aligned}
&\left(f,u_N\right) \\
=& \int_{0}^{1}\int_{0}^{1} u_N(x_1, x_2)\Bigg(- x_2(x_2 - 1)
\left[ \left(x_1 - \frac{1}{2}\right)^{-\frac{2}{3}}
\left(\frac{70}{9}x_1\left(x_1 - 1\right) + \frac{11}{6}\right) 
+ 2\left(x_2 - \frac{1}{2}\right)^{\frac{4}{3}}\right] \\
&- x_1(x_1 - 1)\left[ \left(x_2 - \frac{1}{2}\right)^{-\frac{2}{3}}
\left(\frac{70}{9}x_2\left(x_2 - 1\right) + \frac{11}{6}\right) 
+ 2\left(x_1 - \frac{1}{2}\right)^{\frac{4}{3}}\right] \Bigg) dx_1dx_2 \\
=& I_1 + I_2 + I_3 + I_4.
\end{aligned}
\]
First, we define
\[
\begin{aligned}
\widehat{u}_1(x_1,x_2) &= - x_2(x_2 - 1)\left(\frac{70}{9} x_1 \left(x_1 - 1 \right) 
+ \frac{11}{6}\right)u_N(x_1, x_2), \\
\widehat{u}_2(x_1,x_2) &= - x_1(x_1 - 1)\left(\frac{70}{9} x_2 \left(x_2 - 1 \right) 
+ \frac{11}{6}\right)u_N(x_1, x_2).
\end{aligned}
\]
Then $I_1$, $I_2$, $I_3$, and $I_4$ are expressed as follows:
\[
\begin{aligned}
I_1 &= \int_{0}^{1} \int_{0}^{1} \left(x_1 - \frac{1}{2}\right)^{-\frac{2}{3}} 
\widehat{u}_1(x_1,x_2)dx_1dx_2, \\
I_2 &= \int_{0}^{1} \int_{0}^{1}\left(x_2 - \frac{1}{2}\right)^{-\frac{2}{3}} 
\widehat{u}_2(x_1,x_2)dx_1dx_2, \\
I_3 &= \int_{0}^{1}\int_{0}^{1}- 2x_2(x_2 - 1)
\left(x_2 - \frac{1}{2}\right)^{\frac{4}{3}}u_N(x_1,x_2)dx_1dx_2, \\
I_4 &= \int_{0}^{1}\int_{0}^{1}- 2x_1(x_1 - 1)
\left(x_1 - \frac{1}{2}\right)^{\frac{4}{3}}u_N(x_1,x_2)dx_1dx_2.
\end{aligned}
\]
Now we introduce the computing method for $I_1$, and $I_2$ can be computed with the 
similar way for $I_1$, while $I_3$ and $I_4$ can be computed 
using the Gauss-Legendre quadrature on $[0,1]^2$. The computation of $I_1$ is as follows:
\[
\begin{aligned}
I_1 &= \int_{0}^{1} \int_{0}^{1} \left(x_1 - \frac{1}{2}\right)^{-\frac{2}{3}}
\widehat{u}_1(x_1,x_2)dx_1dx_2 \\
&=\int_{0}^{1}\int_{0}^{\frac{1}{2}}\left(x_1 - \frac{1}{2}\right)^{-\frac{2}{3}} 
\widehat{u}_1(x_1,x_2)dx_1dx_2 + \int_{0}^{1}\int_{\frac{1}{2}}^{1}
\left(x_1 - \frac{1}{2}\right)^{-\frac{2}{3}} \widehat{u}_1(x_1,x_2)dx_1dx_2 \\
&= 4^{-\frac{1}{3}}\int_{0}^{1} \int_{-1}^{1}(1 - x_1)^{-\frac{2}{3}} 
\widehat{u}_1\left(\frac{x_1 + 1}{4}, x_2\right)dx_1dx_2 \\
&\ \ + 4^{-\frac{1}{3}}\int_{0}^{1} \int_{-1}^{1}(1 + x_1)^{-\frac{2}{3}} 
\widehat{u}_1\left(\frac{x_1 + 3}{4}, x_2\right)dx_1dx_2.
\end{aligned}
\]
For the first integral term, we treat $(1 - x_1)^{-\frac{2}{3}}$ as 
the weight function with $(\alpha, \beta) = (-\frac{2}{3},0)$. 
Therefore, the integration method is as follows: in the $x_1$ direction, 
we use the $(-\frac{2}{3},0)$-Gauss-Jacobi quadrature, and in the $x_2$ direction, 
we use the composite Gauss-Legendre quadrature over $[0,1]$ to compute 
the two-dimensional integral of the function $\widehat{u}_1\left(\frac{x_1 + 1}{4}, x_2\right)$. 
For the second integral term, we treat $(1 + x_1)^{-\frac{2}{3}}$ as the weight 
function with $(\alpha, \beta) = (0,-\frac{2}{3})$. Therefore, the integration 
method is as follows: in the $x_1$ direction, 
we use the $(0,-\frac{2}{3})$-Gauss-Jacobi quadrature, and in the $x_2$ 
direction, we use the composite Gauss-Legendre quadrature over $[0,1]$ 
to compute the two-dimensional integral of the 
function $\widehat{u}_1\left(\frac{x_1 + 3}{4}, x_2\right)$.

We have introduced the integration methods used for various parts of the loss function. 
Next, we describe the numerical experiment settings for this example. 

First, we select a TNN structure with a rank parameter $p=100$. Each sub-network is 
a fully-connected neural network (FNN) with three hidden layers, 
each of width 100. When computing the integrals in the loss function, 
we use 200 quadrature points for the Gauss-Jacobi quadrature method. 
For the Gauss-Legendre quadrature method, the integration region 
is divided into 100 subintervals, with 8 Gauss-Legendre quadrature points 
selected in each subinterval.

\begin{figure}[htb!]
\centering
\includegraphics[width=15cm,height=4cm]{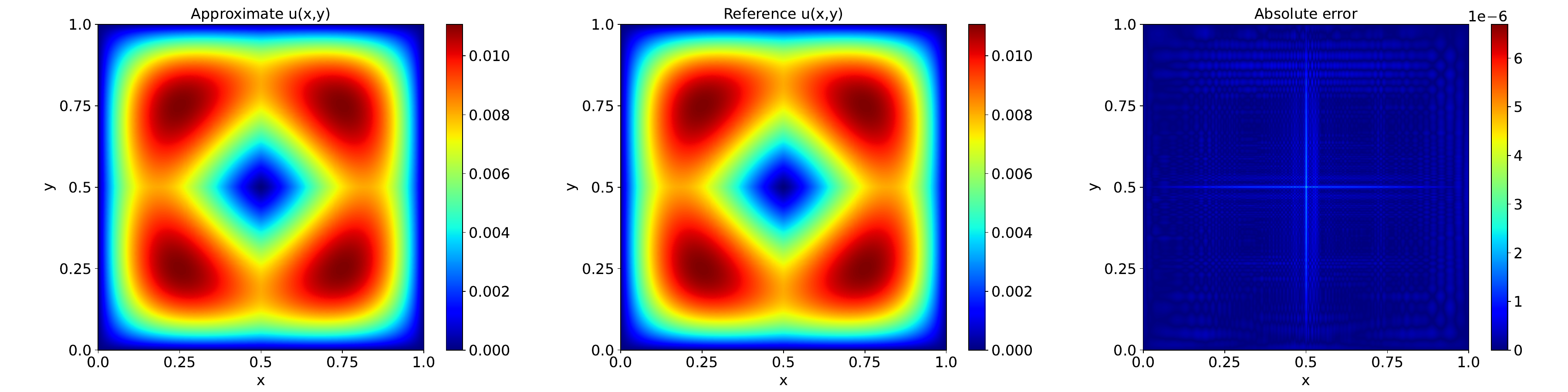}
\caption{Numerical results for solving problem (\ref{Laplace_Singular}): 
the left panel shows the approximate solution, the middle panel shows 
the exact solution, and the right panel shows the error distribution.}
\label{fig_example_1_Ritz_Singular_abs}
\end{figure}

The Adam optimizer with a learning rate of 0.001 is used for the first 50,000 training steps, 
followed by the LBFGS optimizer with a learning rate of 0.1 for the following 
10,000 steps to obtain the final result. 
Figure \ref{fig_example_1_Ritz_Singular_abs} shows the comparison between 
the numerical solution obtained by the algorithm and the exact solution, 
as well as the error distribution. The left panel shows the approximate solution, 
the middle panel shows the exact solution, and the right panel 
shows the error distribution. The $L^2$ relative error between 
the numerical solution and the exact solution is $2.893995 \times 10^{-5}$.

Additionally, we calculate the relative error of the numerical solution $u_N$ against 
the exact solution $u^*$ on a uniformly distributed $1001 \times 1001$ 
test points $\{x^{(k)}, k=1,\dots,K\}$ over $[0,1]^2$. 
The relative error at the test points is defined as
\[
\frac{\sqrt{\mathop\sum\limits_{k=1}^{K}\left (u_N(x^{(k)}) 
- u^*(x^{(k)})\right )^2}}{\sqrt{\mathop\sum\limits_{k=1}^{K}
\left (u^*(x^{(k)})\right )^2}}, 
\]
which yields $2.190126 \times 10^{-5}$.

We also solved this problem using the h-adpative finite element method (FEM) 
with linear and quadratic element. Using adaptive refinement technique, 
we abtain the finite element space of 128888 elements for linear element 
and 86926 elements for quadratic element, 
achieving the accuracy of $2.320817 \times 10^{-3}$ 
and $2.333962 \times 10^{-3}$ for the finite element 
solution with respect to the relative error in the $L^2$ norm.
This demonstrates that the machine learning method can achieve 
comparable accuracy to the finite element method.

To further verify the assertion that the integration accuracy 
affects the accuracy of the machine learning solution, 
we used the same TNN and optimization process as in the previous 
numerical example. However, we use 4000 Gauss-Legendre quadrature points 
per dimension to compute the loss function without specifically addressing 
the singularity in the loss function. 
The resulting relative error on the uniformly distributed test 
points is $4.601047 \times 10^{-2}$. 

We also refined the quadrature points near the singularity regions 
[0.49, 0.51] and [0.499, 0.501], achieving relative errors of $1.249979 \times 10^{-2}$ 
and $5.822244 \times 10^{-3}$, respectively. This shows that as the 
integration accuracy increases, the algorithm's precision also improves. 
The introduction of the Gauss-Jacobi quadrature method effectively handles 
the singularities in the integrand, resulting in the highest accuracy among 
these tests.

We also do the numerical experiment by using ResNet as the sub-network 
of the TNN to enhance training stability. We selected a rank parameter $p=100$ 
and constructed ResNet sub-networks with 9 hidden layers, each of width 50, 
with residual connections between every two hidden layers. 
The quadrature settings are chosen as the same as that for FNN case. 
The Adam optimizer with a learning rate of 0.003 was used for the first 10,000 
training steps, followed by 1000 steps of the LBFGS optimizer with a learning rate of 1  
to obtain the final result. The relative $L^2$ error between the numerical solution 
and the exact solution is $1.797906 \times 10^{-5}$, and the relative error on the 
test points is $7.295925 \times 10^{-6}$. 

\begin{figure}[htb!]
\centering
\includegraphics[width=6cm,height=6cm]{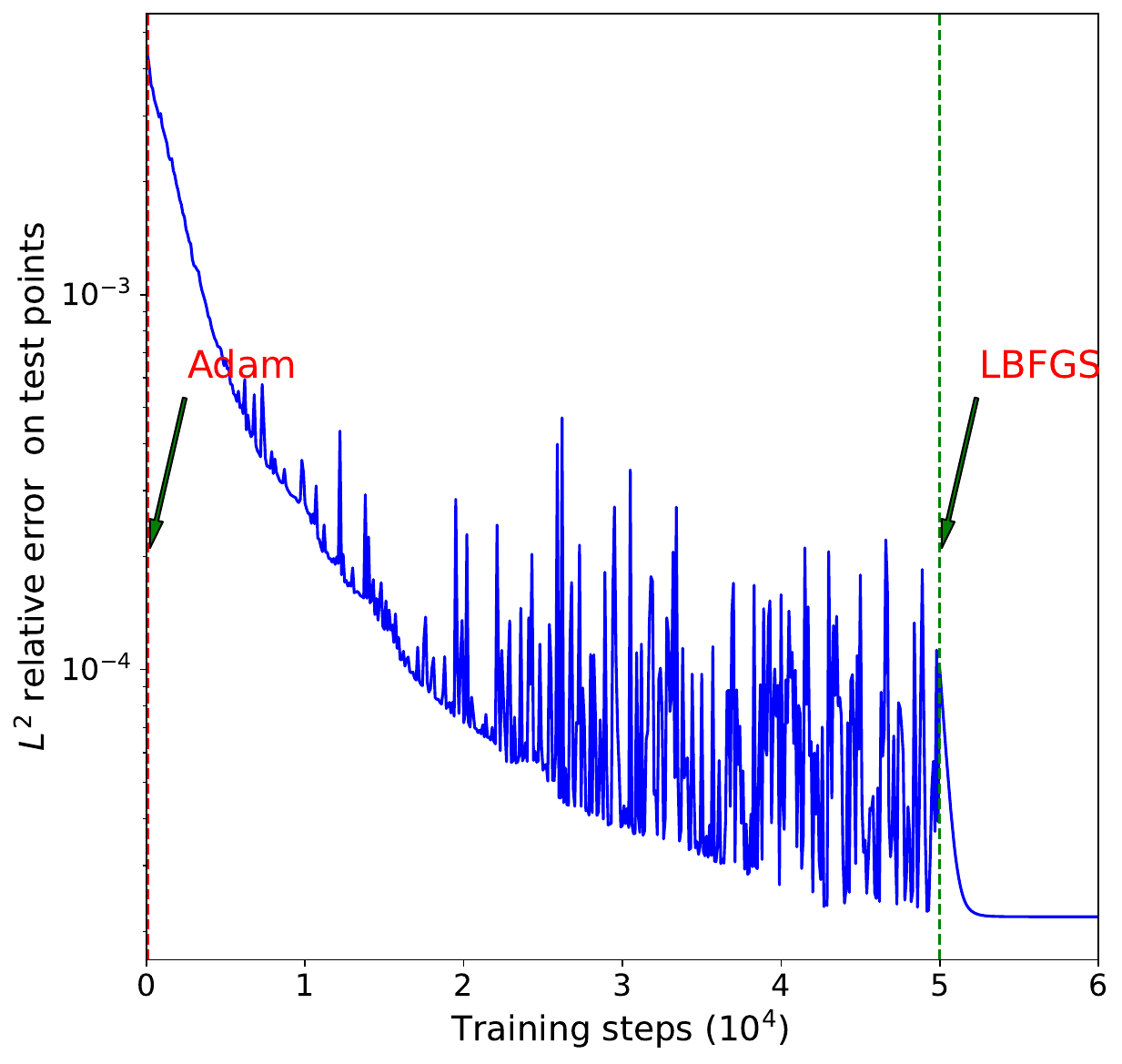}
\includegraphics[width=6cm,height=6cm]{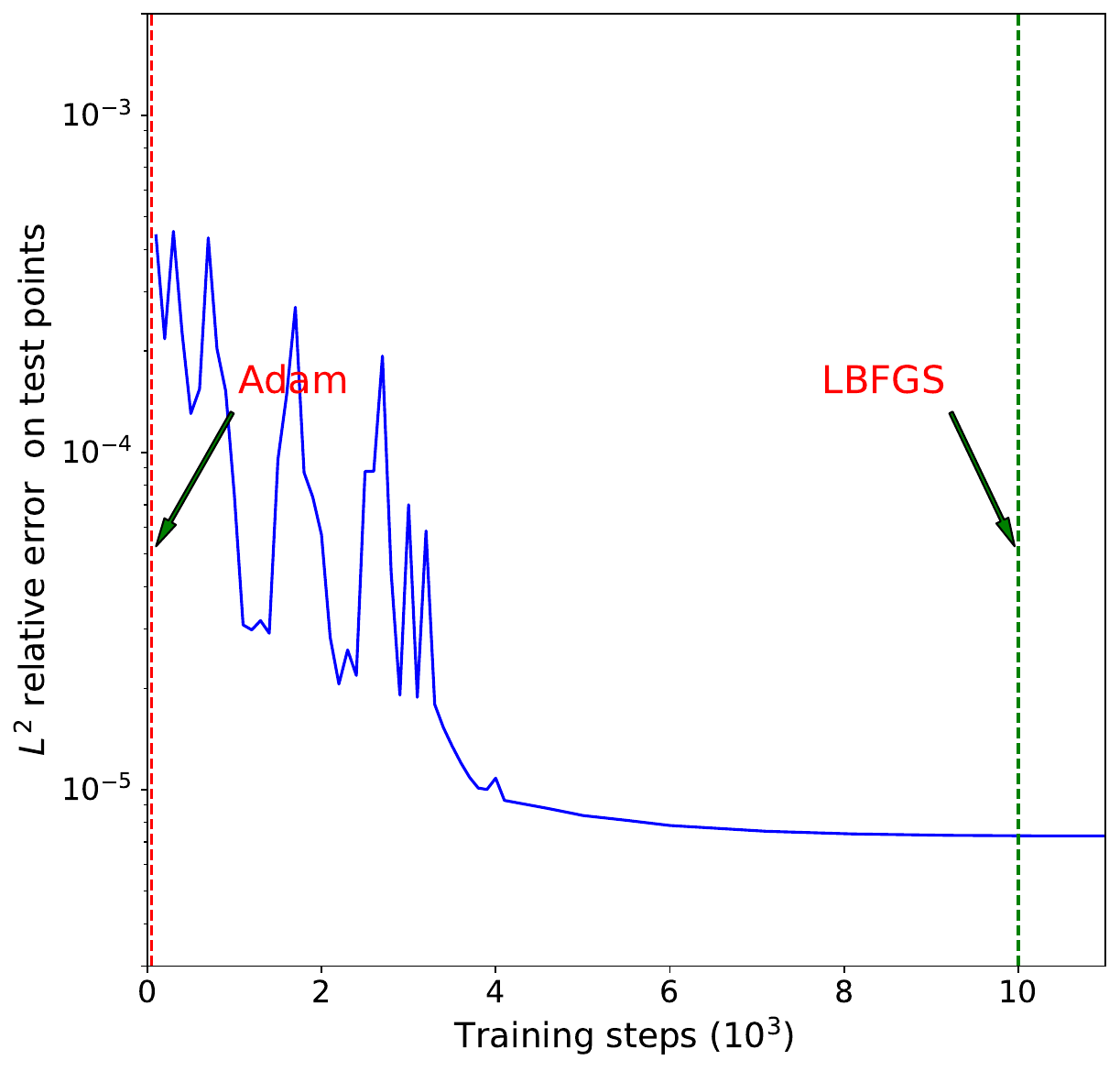}
\caption{The relative error on test points during training for 
TNNs with different sub-network structures. 
The left panel shows FNN, and the right panel shows ResNet.}
\label{fig_example_1_Ritz_Singular_relative}
\end{figure}

Figure \ref{fig_example_1_Ritz_Singular_relative} shows the relative error 
on the test points during training for TNNs with FNN and ResNet sub-networks. 
It can be seen that the training process of ResNet is more stable 
than that of FNN, which aligns with the current understanding in the research community.

Finally, we also investigate the performance of the subspace which is built 
by a two-dimensional FNN with 4 hidden layers, each of width 64, 
and a rank parameter of $p=100$. The two-dimensional quadrature points 
were generated as the Cartesian product of one-dimensional quadrature points. 
The number of quadrature points, optimization process, 
and optimization parameters were consistent with that for the TNN case. 
The resulting relative $L^2$ error between the numerical solution and 
the exact solution is $4.567496 \times 10^{-5}$, and the relative error on 
the test points is $4.621387 \times 10^{-5}$.

\section{Elliptic problems with singularity or discontinuous coefficient}\label{Section_Interface}
In this section, we consider solving a second-order elliptic equation with discontinuous 
coefficients using the machine learning method. 
The primary difference from the previous section lies in the fact that the gradient 
of the solution approximated by the neural network (NN) is not globally continuous. 
Instead, it exhibits discontinuities across interfaces and is a piecewise continuous 
divergence function. This necessitates corresponding modifications when designing 
the a posterior error estimator for the neural network approximation.
In this section, we also provide several numerical examples to verify 
the efficiency and accuracy of the machine learning method based 
on the adaptive neural network subspace proposed in this section.

Interface problems arise in various physical and engineering applications, 
such as materials science, fluid dynamics, and electromagnetics, 
where the coefficients of partial differential equations (PDEs) 
are discontinuous across material interfaces that separate different 
physical regions. When the interface is sufficiently smooth, the solution 
of the interface problem is highly smooth within the subdomains 
occupied by individual materials. However, its global smoothness is typically 
very low. As a result, numerical methods for interface problems have 
become a popular topic in the field of scientific and engineering computation. 
Due to the low global smoothness and irregular geometry 
of the interface, achieving high accuracy by the finite element method 
is often challenging. Various finite element 
methods \cite{chen2023arbitrarily,chen2024arbitrarily,chen1998finite} 
and machine learning methods based on neural 
networks \cite{wu2022inn,xie2024physics,yao2023deep} have been 
developed for solving interface problems. 


This subsection will demonstrate how to design the neural network 
based machine learning method for solving interface problems. 
Additionally, we will also investigate the performances for both Ritz-type 
loss function and the one by the a posterior error estimator. 

\begin{figure}[ht!]
\centering
\begin{tikzpicture}
\draw[help lines,color=gray!160,step=10pt,very thick, xshift =100pt] (4,4) rectangle(8,8);
\draw[help lines,color=gray!160,step=10pt,very thick, xshift =100pt] (6,6) circle(1); 
\node at (9.2,6) [right] {$\Omega_1$}; 
\node at (7.8,4.5) [right] {$\Omega_2$}; 
\draw[help lines,color=gray!160,step=10pt,very thick, xshift =100pt,->] (5.3,6.7) -- (4.75,7.25);
\node at (6.8,6.5) [right] {$\partial\Omega$}; 
\node at (10.5,6.25) [right] {$\Gamma$}; 
\node at (8,6.85) [right] {$\mathbf n$}; 
\end{tikzpicture}
\caption{The boundary $\partial\Omega$, subdomains $\Omega_1$, $\Omega_2$ and interface $\Gamma$}\label{Figure_Interface}
\end{figure}
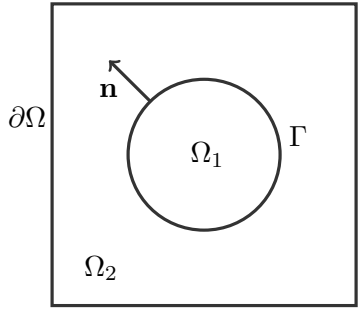

In order to show the way to build the neural network based machine learning method with 
the a posteriori error estimate as the loss functions, we consider solving the interface problem.  
Let $\Omega\subset \mathbb R^2$ be a convex polygonal and $\Omega_1\subset \Omega$ be an open domain 
with Lipschitz boundary $\Gamma=\partial\Omega_1\subset \Omega$. 
Let $\Omega_2=\Omega\backslash \Omega_1$ (see Figure \ref{Figure_Interface}). 

Here, we consider the following elliptic interface problem: Find $u\in H_0^1(\Omega)$ such that 
\begin{eqnarray}\label{Interface_Problem}
\left\{
\begin{array}{rcl}
-\nabla\cdot(\alpha \nabla u) +\beta u&=& f,\ \ \ {\rm in}\ \Omega_1\cup \Omega_2,\\
u&=&0,\ \ \ {\rm on}\ \partial\Omega,\\
\left[ u \right]=0,\ \ \left[\alpha\frac{\partial u}{\partial\bf n}\right]&=&g,\ \ \ {\rm across}\ \Gamma, 
\end{array}
\right.
\end{eqnarray}
where $[v]$ is the jump of a quantity $v$ across the interface $\Gamma$ and $\mathbf n$ the unit outward normal 
to the boundary $\partial\Omega_1$. For definiteness, let $[v](x)=v_1(x)-v_2(x)$, $x\in \Gamma$, with $v_1$ 
and $v_2$ the restrictions of $v$ on $\Omega_1$ and $\Omega_2$, respectively. 
For the ease of description, we assume that the coefficient
function $\beta\geq 0$ and  $\alpha$ is positive and piecewise constant, i.e.  
$$\alpha(x) = \alpha_1\ \ {\rm for}\  x \in \Omega_1;\ \ 
 \alpha(x) = \alpha_2\ \ {\rm for}\  x\in\Omega_2.$$
But the results of this section can be easily extended to more general elliptic
interface problems. 

For building the trial function, we use the functions $\omega(x)$ and $\gamma(x)$ 
to denote the boundaries $\partial\Omega$ and $\Gamma$ implicity, i.e. 
\begin{eqnarray*}
\partial\Omega:=\left\{x\in \mathbb R^2| \omega(x)=0\right\},\ \ \ \ 
\Gamma:=\left\{x\in \mathbb R^2| \gamma(x)=0\right\}. 
\end{eqnarray*}

In order to derive the a posteriori error estimates as the loss function, 
we first define the 
corresponding variational form for (\ref{Interface_Problem}) as follows: 
Find $u\in H_0^1(\Omega)$ such that 
\begin{eqnarray}\label{Interface_Weak}
 a(u,v)=(f,v)+\langle g,v \rangle,\ \ \ \forall v\in H_0^1(\Omega),
\end{eqnarray}
where the bilinear forms $a(\cdot,\cdot)$, $(\cdot,\cdot)$ 
and $\langle\cdot,\cdot\rangle$ are defined as follows
\begin{eqnarray*}
a(u,v)=\int_\Omega \left( \nabla v \cdot \alpha \nabla u+\beta uv \right)d\Omega,\ \ \ 
(f,v)=\int_\Omega fvd\Omega,\ \ \ \langle g,v\rangle =\int_\Gamma gvds.
\end{eqnarray*}
Based on the bilinear form $a(\cdot,\cdot)$, 
we can define the energy norm $\|\cdot\|_a$ as follows
\begin{eqnarray*}
\|v\|_a :=\sqrt{a(v,v)},\ \ \ \forall v\in H_0^1(\Omega).
\end{eqnarray*}

It is easy to know the Ritz form for the problem (\ref{Interface_Weak}) is 
\begin{eqnarray}\label{ritz_loss_interface}
u = \inf_{u\in V}\left(\frac{1}{2}a(v,v)-(f,v)-\langle g,v\rangle\right). 
\end{eqnarray}
Then the Ritz-type loss function for the NN function $\psi$ can be defined as follows
\begin{eqnarray}
\mathcal L(\psi) = \frac{1}{2}a(\psi,\psi)-(f,\psi)-\langle g,\psi\rangle.
\end{eqnarray}

Similarly, we can also derive the loss function by using the  a posteriori 
error estimates for the interface problem. 
It is easy to know $\nabla u\in H({\rm div},\Omega_1)\cap H({\rm div}, \Omega_2)$ 
but $\nabla u\not\in H({\rm div},\Omega)$ \cite{chen1998finite}. 
Before deriving the a posteriori error estimate, we introduce the 
following formula of integration by parts
\begin{eqnarray}\label{Integration_Part}
(\mathbf p, \nabla v) +({\rm div}\mathbf p, v)=\langle [\mathbf n\cdot\mathbf p],  v\rangle,\ \ \ 
\forall \mathbf p\in H({\rm div},\Omega_1)\cap H({\rm div},\Omega_2),\ \ \forall v\in H_0^1(\Omega).
\end{eqnarray}

In order to measure the relation between two norms $\|\cdot\|_{0,\Gamma}$ and $\|\cdot\|_a$, 
we define a Steklov eigenvalue problem: Find $(\lambda_\Gamma, u)\in \mathbb R\times H_0^1(\Omega)$ such that 
\begin{eqnarray}\label{Steklov_Problem}
\left\{
\begin{array}{rcl}
-\nabla\cdot(\alpha \nabla u) +\beta u&=& 0,\ \ \ {\rm in}\ \Omega_1\cup \Omega_2,\\
u&=&0,\ \ \ {\rm on}\ \partial\Omega,\\
\left[ u \right]=0,\ \ \left[\alpha\frac{\partial u}{\partial\bf n}\right]&=&\lambda_\Gamma u,\ \ \ {\rm across}\ \Gamma.
\end{array}
\right.
\end{eqnarray}
The corresponding weak form is defined as follows: Find $(\lambda_\Gamma, u)\in \mathbb R\times H_0^1(\Omega)$ such that 
\begin{eqnarray}\label{Steklov_Weak}
a(u,v)=\lambda_\Gamma \langle u,v\rangle,\ \ \ \forall v\in H_0^1(\Omega).
\end{eqnarray}
It is well-known that eigenvalue problem (\ref{Steklov_Weak}) has the eigenvalue series 
$$0< \lambda_{1,\Gamma} \leq \lambda_{2,\Gamma}\leq \cdots\lambda_{km\Gamma}\leq \cdots\rightarrow +\infty.$$
Then we have the following inequality 
\begin{eqnarray}\label{Steklov_Inequality}
\frac{\|v\|_a^2}{\|v\|_{0,\Gamma}^2}\geq \lambda_{1,\Gamma}\ \ \ \forall 0\neq v\in H_0^1(\Omega). 
\end{eqnarray}

\begin{theorem}
Assume $u\in V$ is the exact solution of (\ref{Interface_Weak}) and the trial function $u_N\in H_0^1(\Omega)$. 
The following upper bound holds
\begin{eqnarray}\label{Upper_Bound_3}
\|u-u_N\|_a&\leq& \sqrt{1+\frac{1}{\lambda_{\rm min}}+\frac{1}{\lambda_{1,\Gamma}}}\eta(u_N, \mathbf p), \ \ \ 
\forall \mathbf p\in \bigcap_{i=1}^2 H({\rm div};\Omega_i),  
\end{eqnarray}
where $\eta(u_p,\mathbf p)$ is defined as follows
\begin{eqnarray}\label{Eta_2}
\eta(u_N, \mathbf p):= \left(\|f-\beta u_N+{\rm div}\mathbf p\|_0^2 +\|\alpha^{-\frac{1}{2}}(\mathbf p-\alpha\nabla u_N)\|_0^2 
+ \|g-[\mathbf n\cdot\mathbf p]\|_{0,\Gamma}^2\right)^{\frac{1}{2}}.
\end{eqnarray}
\end{theorem}
\begin{proof}
Combining (\ref{Interface_Weak}), (\ref{Integration_Part}) and (\ref{Steklov_Inequality}), 
for any function $u_N\in H_0^1(\Omega)$, we can define $w = u-u_N \in V$, and the following inequalities hold for 
any $\mathbf p\in H({\rm div},\Omega_1)\cap H({\rm div},\Omega_2)$
\begin{eqnarray}\label{Hyper_Circle}
&&\|u-u_N\|_a^2 = a(u-u_N,w) = (f,w)+\langle g,w\rangle -a(u_N,w)\nonumber\\
&&=(f,w)+\langle g,w\rangle - (\alpha \nabla u_N, \nabla w)-(\beta u_n,w) + (\mathbf p, \nabla w) +({\rm div}\mathbf p, w)
-\langle [\mathbf n\cdot\mathbf p],  w\rangle\nonumber\\
&&=(f-\beta u_N+{\rm div}\mathbf p,w)+(\mathbf p-\alpha\nabla u_N,\nabla w)+ \langle g- [\mathbf n\cdot\mathbf p],  w\rangle\nonumber\\
&&\leq \|f-\beta u_N+{\rm div}\mathbf p\|_0 \|w\|_0 
+ \|\alpha^{-\frac{1}{2}}(\mathbf p-\alpha\nabla u_N)\|_0\|\alpha^{\frac{1}{2}}\nabla w\|_0\nonumber\\
&&\ \ \ \  + \|g- [\mathbf n\cdot\mathbf p]\|_{0,\Gamma}\|w\|_{0,\Gamma}\nonumber\\
&&\leq \left(\|f-\beta u_N+{\rm div}\mathbf p\|_0^2+\|\alpha^{-\frac{1}{2}}(\mathbf p-\alpha\nabla u_N)\|_0^2 
+ \|g- [\mathbf n\cdot\mathbf p]\|_{0,\Gamma}^2\right)^{\frac{1}{2}}\nonumber\\
&&\ \ \ \ \times
\left(\|w\|_0^2+\|\alpha^{\frac{1}{2}}\nabla w\|_0^2 + \|w\|_{0,\Gamma}^2\right)^{\frac{1}{2}}\nonumber\\
&&\leq \left(\|f-\beta u_N+{\rm div}\mathbf p\|_0^2+\|\alpha^{-\frac{1}{2}}(\mathbf p-\alpha\nabla u_N)\|_0^2 
+ \|g- [\mathbf n\cdot\mathbf p]\|_{0,\Gamma}^2\right)^{\frac{1}{2}} \nonumber\\
&&\ \ \ \ \times\left(1+\frac{1}{\lambda_{\rm min}}+\frac{1}{\lambda_{1,\Gamma}}\right)^{\frac{1}{2}}\|w\|_a,
\end{eqnarray}
where $\lambda_{\rm min}$ is the constant in the inequality (\ref{Hilbert_Inequality}), 
$\lambda_{1,\Gamma}$ is the minimum eigenvalue of (\ref{Steklov_Weak}). 

Then based on (\ref{Hyper_Circle}), the following inequality holds for any 
$\mathbf p\in \cap_{i=1}^n H({\rm div},\Omega_i)$
\begin{eqnarray*}
\|u-u_N\|_a^2 &\leq& \left(1+\frac{1}{\lambda_{\rm min}}+\frac{1}{\lambda_{1,\Gamma}}\right) 
\Big(\|f-\beta u_N+{\rm div}\mathbf p\|_0^2
+\|\alpha^{-\frac{1}{2}}(\mathbf p-\alpha\nabla u_N)\|_0^2 \nonumber\\
&&\ \ \ \ \ + \|g- [\mathbf n\cdot\mathbf p]\|_{0,\Gamma}^2\Big).
\end{eqnarray*}
This is the desired result (\ref{Upper_Bound_3}) and the proof is complete.
\end{proof}
In order to define the machine learning method, we define the trial function $u_N$ by the neural network functions 
$\Psi(x,\theta)$ and $\Psi_1(x,\theta)$ with the following way 
\begin{eqnarray}\label{Interface_NN}
u_N = \omega(x)\Psi(x,\theta) + \gamma(x)\Psi_1(x,\theta), 
\end{eqnarray}
where $\Psi(x,\theta)$ and $\Psi_1(x, \theta)$ are two neural network functions 
which are defined on $\Omega$ and $\Omega_1$, respectively.  
Then it is easy to know that $[u_N]=0$ and $u_N\in H_0^1(\Omega)$. 
But $\alpha\nabla u_N$ is discontinuous across the interface $\Gamma$ since the supporting 
set of $\Psi_1$ is only $\Omega_1$. Then we only have 
$\alpha \nabla u_N\in H({\rm div},\Omega_1)\cap H({\rm div},\Omega_2)$.

With the help of $\alpha \nabla u_N\in H({\rm div},\Omega_1)\cap H({\rm div},\Omega_2)$, 
we can choose $\mathbf p = \alpha \nabla u_N$ and the following a posteriori error estimate holds
\begin{eqnarray}\label{Posteriori_Error}
\|u-u_N\|_a^2 \leq  \left(1+\frac{1}{\lambda_{\rm min}}+\frac{1}{\lambda_{1,\Gamma}}\right)
\left(\|f-\beta u_N+\nabla\cdot(\alpha \nabla u_N)\|_0^2 
+ \|g- [\mathbf n\cdot(\alpha \nabla u_N)]\|_{0,\Gamma}^2\right).\ \ \ \ 
\end{eqnarray}
Based on the a posteriori error estimate, the loss function can be defined as 
\begin{eqnarray}\label{Interface_Loss}
\mathcal L(u_N) := \eta(u_N,\alpha\nabla u_N)=\|f-\beta u_N+\nabla\cdot(\alpha \nabla u_N)\|_0^2 
+ \|g- [\mathbf n\cdot(\alpha \nabla u_N)]\|_{0,\Gamma}^2.
\end{eqnarray}
Then we can design NN-based machine learning method by using 
Algorithm \ref{Algorithm_Subspace} with the loss function defined by (\ref{Interface_Loss}) 
for solving interface problem (\ref{Interface_Problem}). 

\subsection{Typical two-material diffusion problems}

Consider the following diffusion problem with two materials in the computational domain:
\begin{equation}\label{interface_example_1}
    \left\{
    \begin{array}{rcll}
    -\nabla \cdot(\alpha\nabla u)&=&f, &{\rm in}\ \Omega_1 \cup \Omega_2,  \\ 
    u&=&0, & {\rm on}\  \partial \Omega, \\
    \left[ u \right]=0,\ \ \left[\alpha\frac{\partial u}{\partial\bf n}\right]&=&g,
    \ \ &{\rm across}\ \Gamma
    \end{array}
    \right.
\end{equation}
where $\Omega=(0,1) \times(0,1)$, the discontinuous coefficient $\alpha$ is defined as 
\begin{eqnarray*}
\alpha(x_1,x_2)=
\left\{
\begin{array}{ll}
\alpha_1, & (x_1, x_2) \in\left(0, \frac{2}{3}\right] \times(0,1)=:\Omega_1, \\
\alpha_2, & (x_1, x_2) \in\left(\frac{2}{3}, 1\right) \times(0,1)=:\Omega_2.
\end{array}
\right.
\end{eqnarray*}
The exact solution of this equation is
\begin{eqnarray}\label{Interface_Exact_1}
u(x_1,x_2)=
\left\{
\begin{array}{ll}
    c_1 \sin   (k_{11}\pi x_1) \sin (k_{12} \pi x_2),  &(x_1, x_2) \in\left(0, \frac{2}{3}\right] \times(0,1), \\
    c_2\sin (k_{21} \pi x_1) \sin (k_{22} \pi x_2),  &(x_1, x_2) \in\left(\frac{2}{3}, 1\right) \times(0,1).
\end{array}
\right.
\end{eqnarray}
The source term $f$ and the flux jump $g$ across the interface can be 
derived from the exact solution.

First, we construct the approximation function $u_N$. Based on the idea of 
the spectral element method, the computational domain is divided into 
two subdomains according to the interface: 
$\Omega_1 := \left(0, \frac{2}{3}\right) \times (0,1)$, 
$\Omega_2 := \left(\frac{2}{3}, 1\right) \times (0,1)$. 
We define three neural networks $u_{N}^{(1)}$, $u_{N}^{(2)}$, $u_{N}^{(3)}$ on the subdomains $\Omega_1$, $\Omega_2$, and the entire domain $\Omega$, respectively. By multiplying these networks with different boundary functions, we obtain $\widehat{u}_{N}^{(1)}$, $\widehat{u}_{N}^{(2)}$, $\widehat{u}_{N}^{(3)}$, which satisfy the Dirichlet boundary condition and the continuity condition across the interface:
\[
\begin{aligned}
\widehat{u}_{N}^{(1)} &= x_1\left(\frac{2}{3} - x_1\right)x_2(1-x_2)u_{N}^{(1)}, \\
\widehat{u}_{N}^{(2)} &= \left(x_1 - \frac{2}{3}\right)(1-x_1)x_2(1-x_2)u_{N}^{(2)}, \\
\widehat{u}_{N}^{(3)} &= x_1(1-x_1)x_2(1-x_2)u_{N}^{(3)}.
\end{aligned}
\]
The final approximation function is given by:
\[
u_N = \widehat{u}_{N}^{(1)} + \widehat{u}_{N}^{(2)} + \widehat{u}_{N}^{(3)}.
\]
It follows that:
\[
u_N =
\begin{cases}
\widehat{u}_{N}^{(1)} + \widehat{u}_{N}^{(3)}, & (x_1, x_2) \in \left(0, \frac{2}{3}\right] \times (0,1), \\
\widehat{u}_{N}^{(2)} + \widehat{u}_{N}^{(3)}, & (x_1, x_2) \in \left(\frac{2}{3}, 1\right) \times (0,1).
\end{cases}
\]

Next, we perform numerical experiments using 
the Ritz type of loss function \eqref{ritz_loss_interface} and the one by the a 
posterior error estimator \eqref{Interface_Loss}.

\subsubsection{Testing the loss functions}\label{Section_test_1}
In this subsection, we investigate the performance of two types of loss functions. 
Specifically, we test the performance of the adaptive neural network subspace 
based machine learning method using the Ritz loss function \eqref{ritz_loss_interface} 
and the a posterior error estimator-based loss function \eqref{Interface_Loss}. 

The parameters for problem \eqref{interface_example_1} are set as follows: 
the true solution parameters are $c_1 = c_2 = 1$, $k_{11} = 1$, $k_{12} = 2$, 
$k_{21} = 4$, and $k_{22} = 2$. 
We test the numerical performance under two different material coefficient settings:
\begin{itemize}
    \item Test (1.1): $\alpha_1 = 4$, $\alpha_2 = 1$;
    \item Test (1.2): $\alpha_1 = 4$, $\alpha_2 = 2$.
\end{itemize}

In Test (1.1), the flux across the interface is continuous, whereas in Test (1.2), 
the flux jump across the interface is non-zero. For both tests, 
we use Tensor Neural Networks (TNNs) 
to construct the neural network space. Each sub-network in the TNN consists 
of a fully connected neural network (FNN) 
with three hidden layers, each of width 50, and employs  $\sin(x)$ as the activation function. 

For integration in subdomains $\Omega_1$ and $\Omega_2$, the one-dimensional regions 
are divided into 100 subintervals, with each subinterval using 16 Gauss-Lobatto quadrature points. 
During training, the Adam optimizer with a learning rate of 0.001 is used for the first 
5000 training steps, the LBFGS optimizer with a learning rate of 0.1 
for the following 100 steps to obtain the final results.

We use the following three metrics to measure the errors 
between the neural network approximation $u_h$ and the exact solution $u^*$:
\begin{enumerate}
\item Relative energy norm error: $e_E := \frac{\|u_h - u^*\|_a}{\|u^*\|_a}$;
\item Relative $L^2$ norm error: $e_{L^2} := \frac{\|u_h - u^*\|_{L^2}}{\|u^*\|_{L^2}}$;
\item Relative error on test points: 
\[
e_{\text{test}} := \frac{\sqrt{\sum\limits_{k=1}^{K}\left (u_h(x^{(k)}) 
- u^*(x^{(k)})\right )^2}}
{\sqrt{\sum\limits_{k=1}^{K}\left (u^*(x^{(k)})\right )^2}},
\]
where the test points $\{x^{(k)}, k=1, \dots, K\}$ are uniformly distributed over $[0,1]^2$ with a resolution of $301 \times 301$ grid points.
\end{enumerate}

Tables \ref{table_interface_1_ritz} and \ref{table_interface_1_post} 
present the results of solving the problem using the two loss functions. 
\begin{table}[htb!]
\begin{center}
\caption{Numerical results using the Ritz-type loss function}\label{table_interface_1_ritz}
\begin{tabular}{cccc}
\hline
&   $e_E$ &  $e_{L^2}$  &  $e_{\text{test}}$  \\
\hline
Test (1) & 1.379726e-06 &   1.514646e-07 &   1.473214e-07 \\
\hline
Test (2) & 3.571384e-07 &  1.642144e-07 &   4.6365045e-08 \\
\hline
\end{tabular}
\end{center}
\end{table}
\begin{table}[htb!]
\begin{center}
\caption{Numerical results using the a posteriori error estimator-based loss 
function}\label{table_interface_1_post}
\begin{tabular}{cccc}
\hline
&   $e_E$ &  $e_{L^2}$  &  $e_{\text{test}}$  \\
\hline
Test (1) & 4.442695e-08 &   1.003098e-08 &   1.331919e-09 \\
\hline
Test (2) & 7.731275e-08 &  1.089644e-07 &   8.876892e-09 \\
\hline
\end{tabular}
\end{center}
\end{table}
From Tables \ref{table_interface_1_ritz} and \ref{table_interface_1_post}, 
it can be observed that the machine learning method based on neural network spaces 
can achieve high accuracy in solving this problem. 
The continuity of the flux across the interface does not significantly 
affect the performance of our algorithm. Moreover, the accuracy for 
using the posterior error estimator-based loss function is slightly higher 
compared to that by using the Ritz loss function.

\begin{figure}[htb]
\centering
\includegraphics[width=15cm,height=4cm]{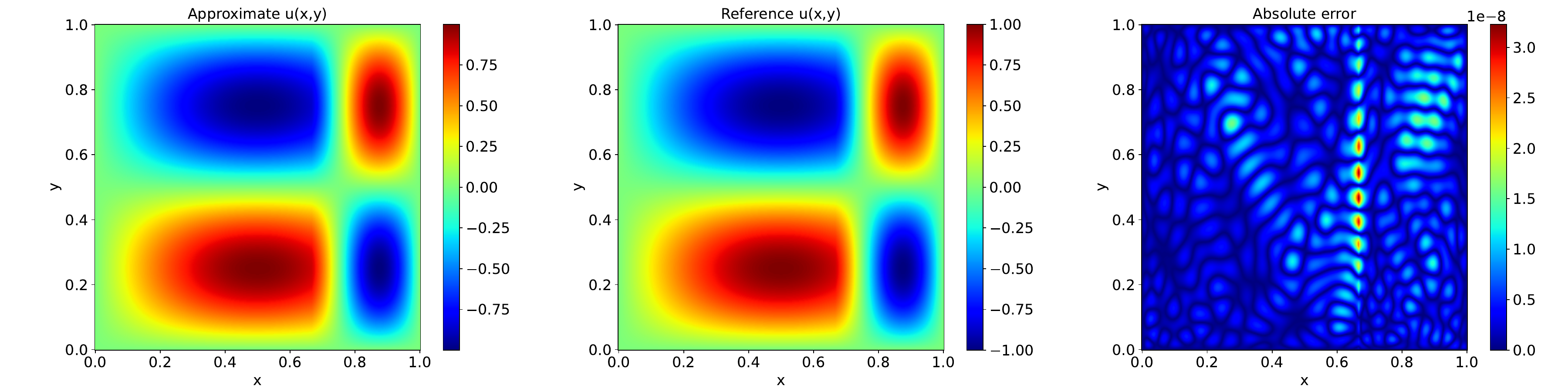}
\caption{Test (1.2): Numerical results of solving (\ref{interface_example_1}) 
using the posterior error estimator-based loss function: 
Left: Approximate solution image, Middle: Exact solution image, 
Right: Error distribution plot}\label{fig_interface_1_2_post_jump}
\end{figure}

\begin{figure}[htb]
\centering
\includegraphics[width=6cm,height=6cm]{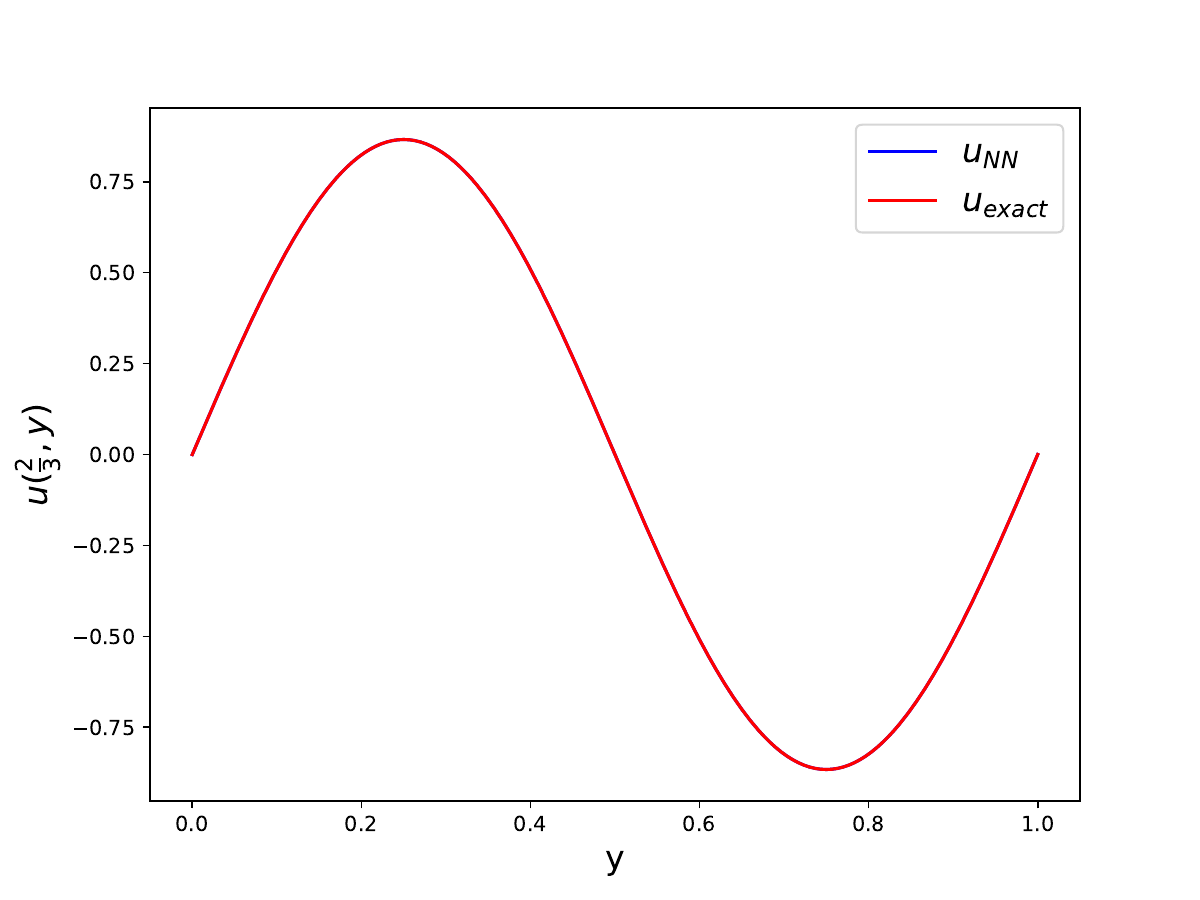}
\includegraphics[width=6cm,height=6cm]{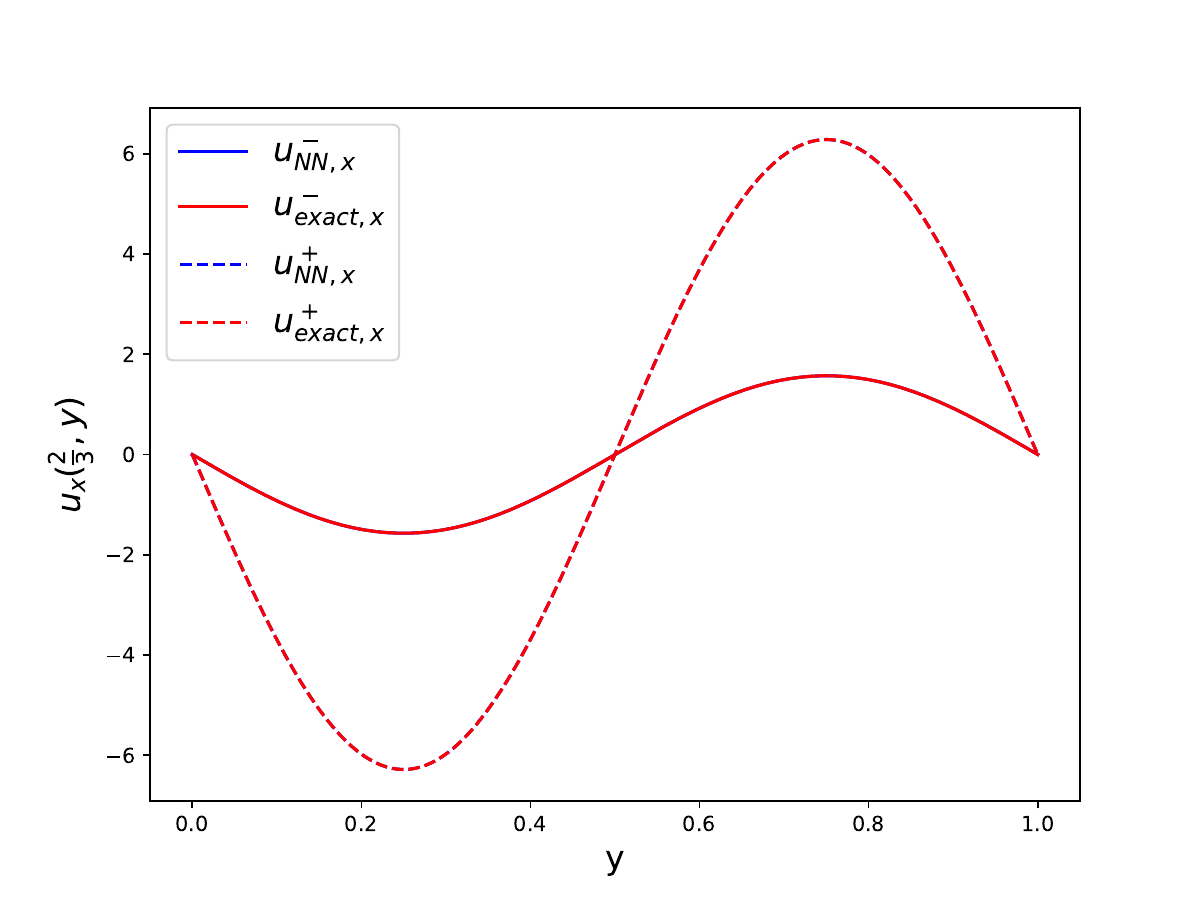}
\caption{Test (1.2): Prediction of the approximate solution at 
the material interface ($x=\frac{2}{3}$) using the posterior 
error estimator-based loss function for (\ref{interface_example_1}), 
Left: Function value prediction, Right: Left and right partial derivative predictions}\label{fig_interface_1_2_post_jump_interface}
\end{figure}

For example, using the posterior error estimator loss function to solve Test (1.2), 
Figure \ref{fig_interface_1_2_post_jump} shows the numerical 
results of the machine learning method designed in this section. 
Figure \ref{fig_interface_1_2_post_jump_interface} displays the 
trajectory of the approximate solution at the material interface ($x=\frac{2}{3}$). 
We see that the machine learning method provides 
accurate approximations for the function values and their partial 
derivatives at the interface.

\subsubsection{Tests with large coefficient differences}\label{Section_test_2}
In this subsection, we will use the machine learning method with the a posterior 
error estimator-based loss function to solve the interface problem 
with large coefficient differences. This is to test the performance 
of the machine learning method when the material 
coefficients have significant differences. 

We set the true solution parameters of the equation (\ref{interface_example_1}) 
as $c_1 = c_2 = 1$, $k_{11} = 1$, $k_{12} = 2$, $k_{21} = 4$, $k_{22} = 2$.  
Here, we test the numerical results for three different material coefficient settings:
\begin{itemize}
\item Test (2.1):  $\alpha_1 = 4 \times 10^{-4}$, $\alpha_2 = 1$;
\item Test (2.2):  $\alpha_1 = 4$, $\alpha_2 = 1$;
\item Test (2.3):  $\alpha_1 = 4000$, $\alpha_2 = 1$.
\end{itemize}

In Test (2.1) and Test (2.3), the solution exhibits flux jumps at the interface, 
and there is a large difference in the coefficients on both sides of 
the interface. In all three tests, we use TNN to construct the neural 
network space. Each TNN has sub-networks selected as FNNs 
with three hidden layers, each layer containing 50 units, and $\sin(x)$ 
as the activation function. When integrating over the subregions 
$\Omega_1$ and $\Omega_2$, we divide the one-dimensional region 
into 100 subintervals, and for each subinterval, we use 16 Gauss-Lobatto 
integration points. During training, we first use the Adam optimizer 
with a learning rate of 0.001 for the first 5000 training steps, 
and then use the LBFGS optimizer with a learning rate of 0.1 for the following 
100 iterations to obtain the final result.

\begin{figure}[htb]
\centering
\includegraphics[width=4.8cm,height=3.84cm]{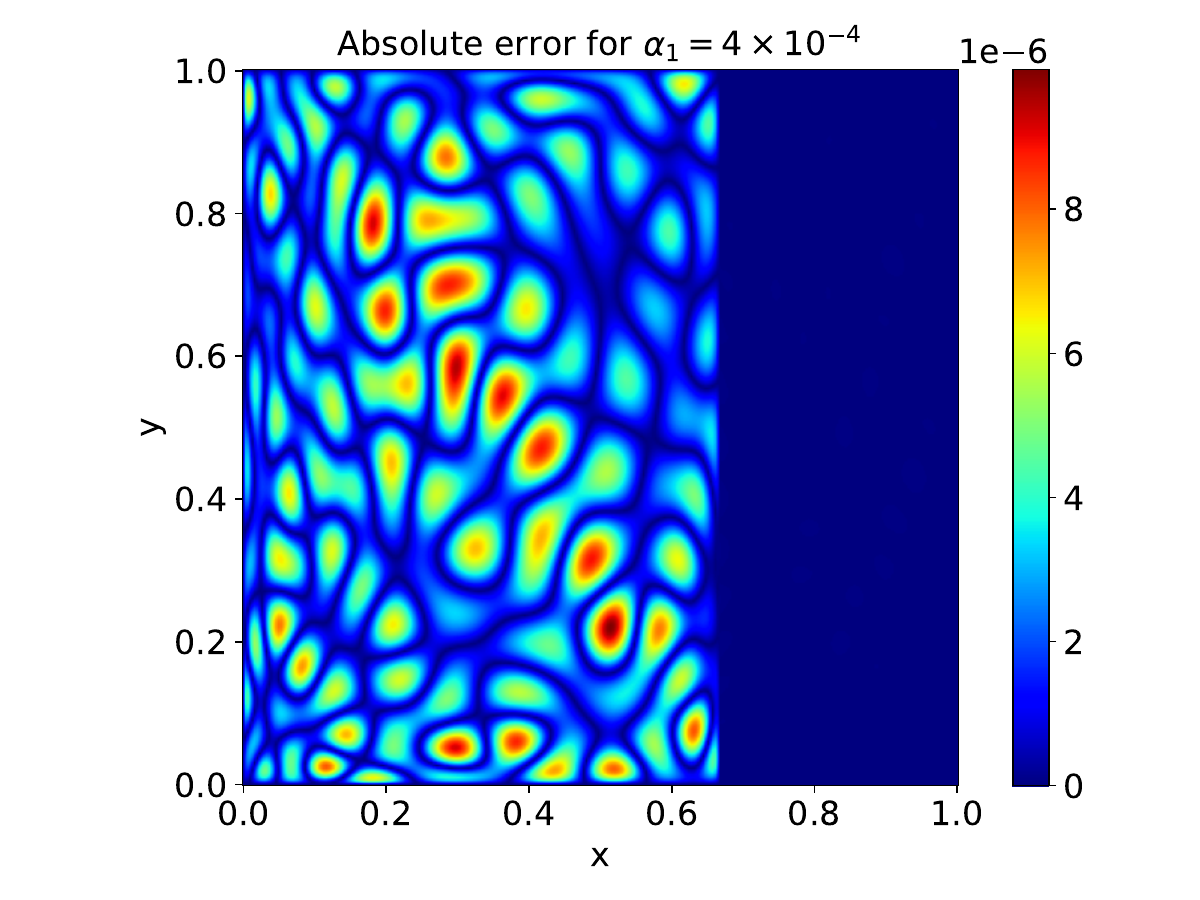}
\includegraphics[width=4.8cm,height=3.84cm]{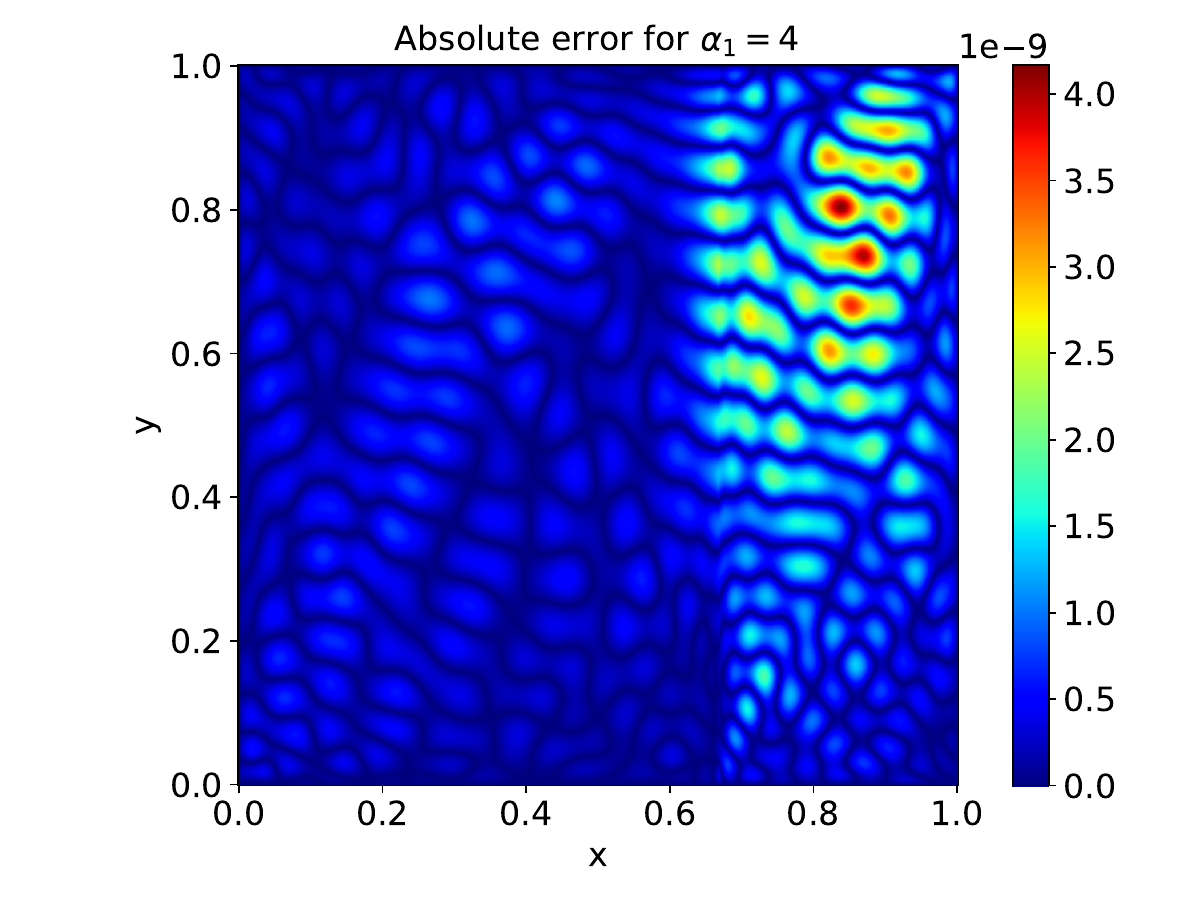}
\includegraphics[width=4.8cm,height=3.84cm]{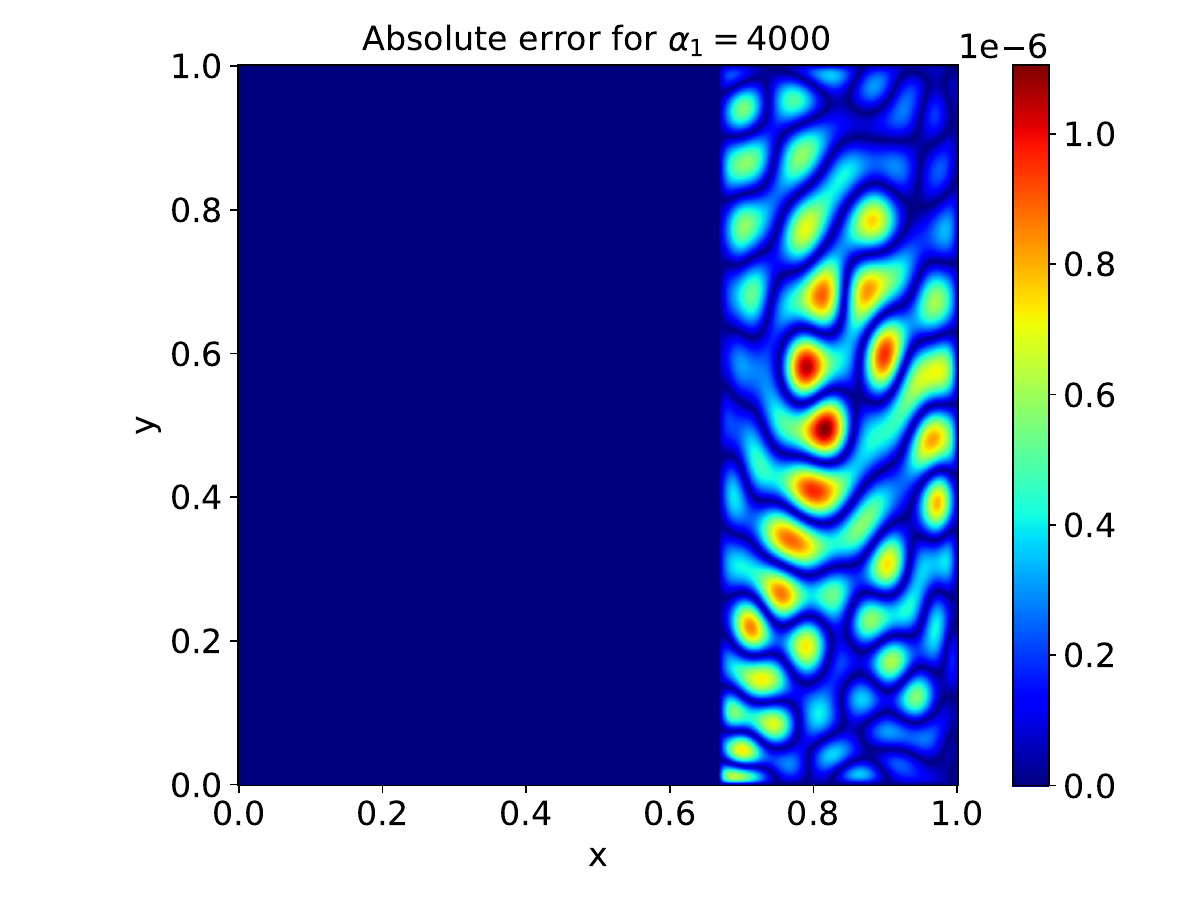}
\caption{Error distribution for the machine learning method  
with three coefficient settings. Left, middle, and right correspond to 
$\alpha_1 = 4 \times 10^{-4}$, $4$, and $4000$, respectively.}\label{fig_interface_2_BD}
\end{figure}

Table \ref{table_interface_2_post_big_difference} presents 
the numerical errors for the three coefficient settings. 
It can be observed that the difficulty of solving the equation 
increases when there is a large coefficient difference, 
but the algorithm still achieves high accuracy. 
Figure \ref{fig_interface_2_BD} shows the error distribution 
for the three coefficient settings. 
From this, we can see that when $\alpha_1 = 4$, the error distribution 
is relatively uniform. However, when the coefficient difference is large, 
in the region with the smaller coefficient, the ellipticity of the equation is poorer, 
leading to larger solution errors, which is consistent with theoretical expectations.
\begin{table}[htb]
\begin{center}
\caption{numerical results for different coefficient differences}
\label{table_interface_2_post_big_difference}
\begin{tabular}{lccc}
\hline 
&   $e_E$ &  $e_{L^2}$  &  $e_{\text{test}}$  \\
\hline
Test (2.1): $\alpha_1 = 4 \times 10^{-4}$   & 7.546526e-07 &  4.964822e-06 & 4.923701e-06 \\
\hline
Test (2.2): $\alpha_1 = 4$                  & 4.442695e-08 &  1.003098e-08 & 1.331919e-09 \\
\hline
Test (2.3): $\alpha_1 = 4000$               & 6.601858e-08 &  3.747886e-07 & 3.731660e-07 \\
\hline
\end{tabular}
\end{center}
\end{table}

\subsubsection{Tests with high-frequency exact solution}\label{Section_test_3}
Then, in this subsection, we test the performance and stability 
of the machine learning method when the exact solution 
contains high-frequency components. For this purpose, 
we consider the following equation: 
\begin{equation}\label{interface_example_12}
\left\{
\begin{array}{rcll}
-\nabla \cdot(\alpha\nabla u) &=& f, & {\rm in} \ \Omega_1\cup\Omega_2,  \\ 
u &=& 0, & {\rm on} \ \partial \Omega, \\
\left[ u \right] = 0, \ \ \left[\alpha \frac{\partial u}{\partial \mathbf{n}}\right] &=& g,
\ \ & {\rm across} \ \Gamma, 
\end{array}
\right.
\end{equation}
where $\Omega = (0,1) \times (0,1)$, and the discontinuous coefficient $\alpha$ is defined as follows
\begin{eqnarray*}
\alpha(x_1,x_2) =
\left\{
\begin{array}{ll}
10, & (x_1, x_2) \in \left(0, \frac{2}{9}\right] \times (0,1)=:\Omega_1, \\
2,  & (x_1, x_2) \in \left(\frac{2}{9}, 1\right) \times (0,1)=:\Omega_2,
\end{array}
\right.
\end{eqnarray*}
The exact solution of the equation is
\begin{eqnarray}\label{Interface_Exact_12}
u(x_1,x_2) =
\left\{
\begin{array}{ll}
 \sin (\pi x_1) \sin (2\pi x_2),  & (x_1, x_2) \in \left(0, \frac{2}{9}\right] \times (0,1), \\
\sin (10 \pi x_1) \sin (2 \pi x_2),  & (x_1, x_2) \in \left(\frac{2}{9}, 1\right) \times (0,1).
\end{array}
\right.
\end{eqnarray}
The corresponding right-hand side $f$ and the flux jump $g$ at the interface 
can be derived from the exact solution. It can be seen from the definition, 
the exact solution has high-frequency components in the right region. 
We label this test as Test (3.1).

We adopt the same approach as in solving (\ref{interface_example_1}), 
defining three neural networks for the two subregions 
and the entire region, and using corresponding boundary factor functions 
to enforce Dirichlet boundary conditions and interface continuity conditions. 
The TNN is used to construct the neural network space, 
selecting FNNs with four hidden layers as sub-networks, 
each with 50 units per layer, and $\sin(x)$ as the activation function. 

The integration method and optimization parameters are exactly the same 
as those used in Test (1.1) in Section \ref{Section_test_1}, 
and we test both loss functions. Table \ref{table_interface_3_error} 
shows the numerical errors obtained using the two loss functions. 
It can be seen that the machine learning method here can effectively 
capture the high-frequency information.

\begin{figure}[htb]
\centering
\includegraphics[width=6cm,height=6cm]{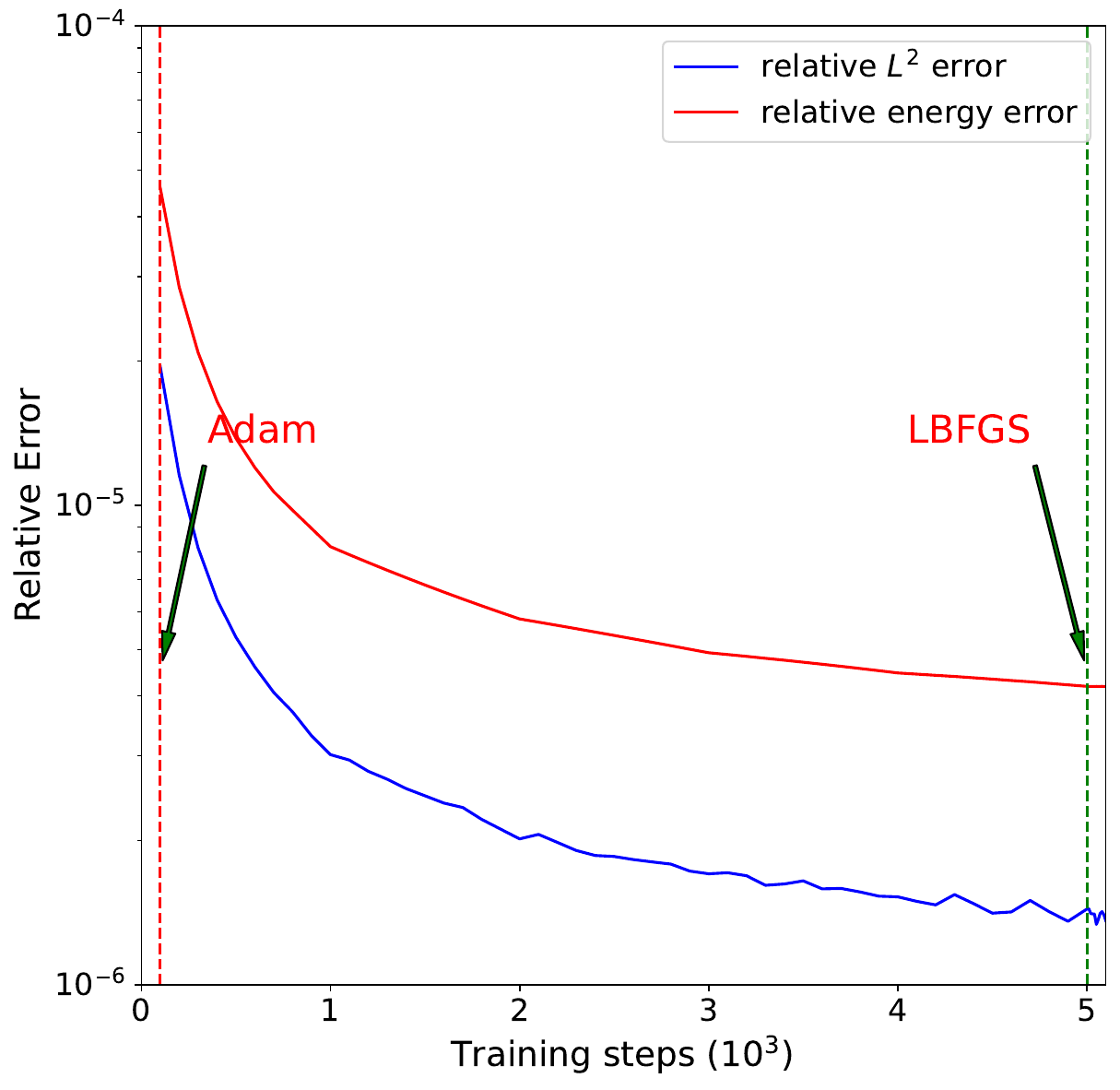}
\includegraphics[width=6cm,height=6cm]{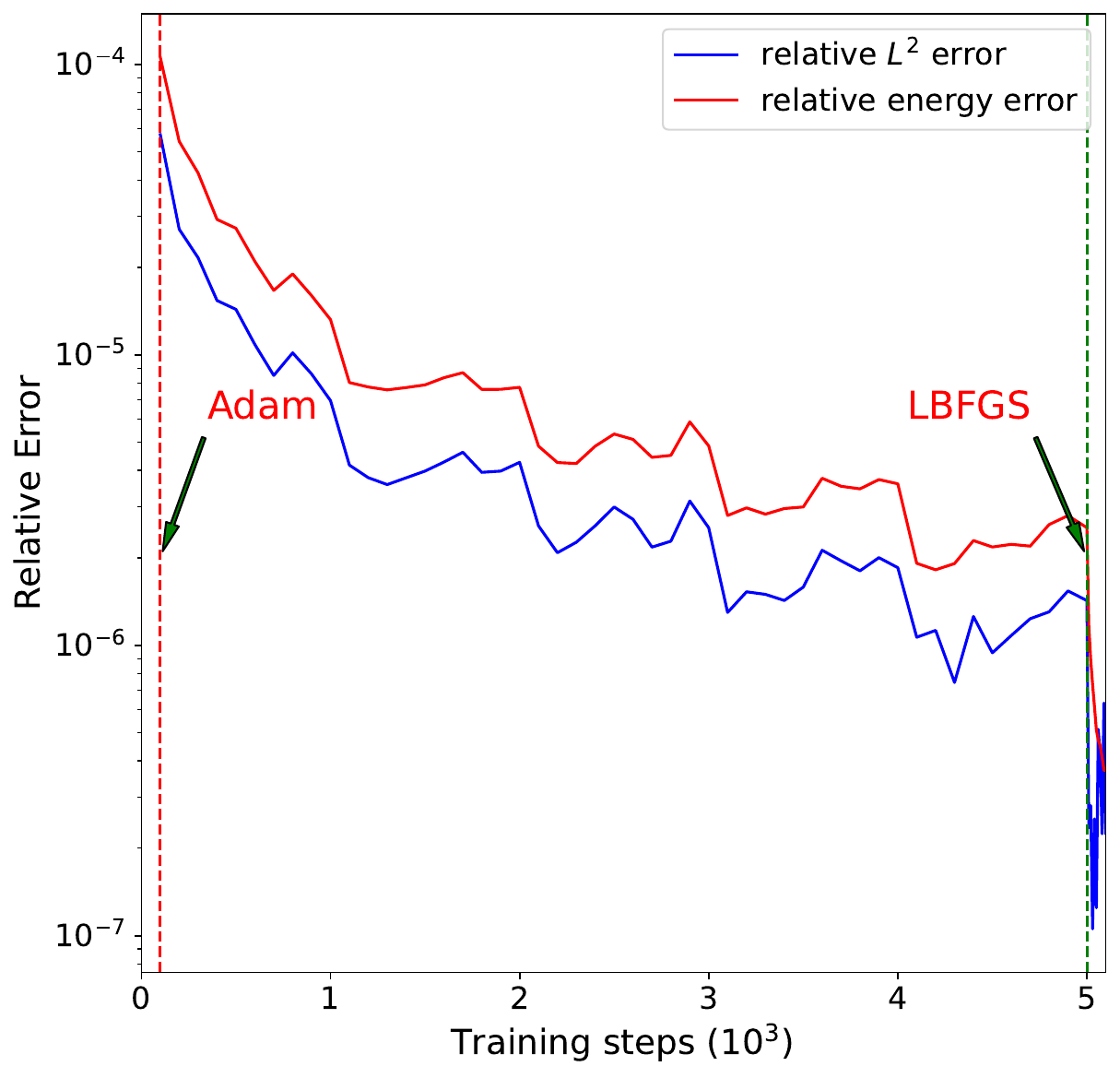}
\caption{Test (3.1): Relative energy error and relative $L^2$ error variation 
during the training process of solving (\ref{interface_example_1}) 
with two loss functions. Left: Ritz-type, right: The a posterior error estimator-type.}\label{fig_interface_3_iteration}
\end{figure}

\begin{figure}[htb]
\centering
\includegraphics[width=15cm,height=4cm]{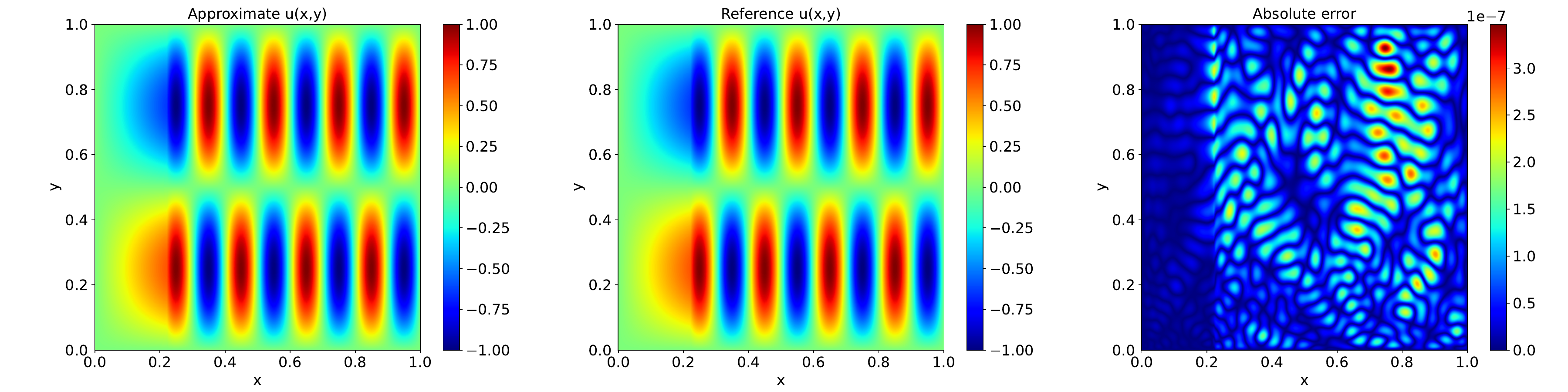}
\caption{Test (3.1): Error distribution obtained by solving (\ref{interface_example_1}) 
with the a posterior error estimator-type loss function.}\label{fig_interface_3_post_jump_abs}
\end{figure}

Figure \ref{fig_interface_3_iteration} shows the relative energy error 
and relative $L^2$ error during the training process. 
Along the steps correspond to the learning process of finding the optimal subspace, 
the machine learning method appears stable. 
Figure \ref{fig_interface_3_post_jump_abs} presents 
the error distribution obtained using the a posterior error estimator-based loss function.
\begin{table}[htb]
\begin{center}
\caption{Numerical results for high-frequency exact solutions using machine learning method}\label{table_interface_3_error}
\begin{tabular}{lccc}
\hline
&   $e_E$ &  $e_{L^2}$  &  $e_{\text{test}}$  \\
\hline
Ritz-type & 4.190128e-06 &   1.359826e-06 &   1.399258e-06 \\
\hline
Posterior error estimator-type:  & 3.702056e-07 &   2.251677e-07 &   1.690076e-07 \\
\hline
\end{tabular}
\end{center}
\end{table}

\subsection{Multi-material diffusion problems}
\begin{figure}[ht]
\centering
\begin{tikzpicture}
\node at (6.6,3.7) [right] {$(-1,-1)$}; 
\node at (10.8,3.7) [right] {$(1,-1)$}; 
\node at (10.8,8.3) [right] {$(1,1)$}; 
\node at (6.7,8.3) [right] {$(-1,1)$}; 
\draw[help lines,color=gray!160,step=10pt,very thick, xshift =100pt,] (4,4) -- (8,4);
\draw[help lines,color=gray!160,step=10pt,very thick, xshift =100pt,] (8,4) -- (8,8);
\draw[help lines,color=gray!160,step=10pt,very thick, xshift =100pt,] (8,8) -- (4,8);
\draw[help lines,color=gray!160,step=10pt,very thick, xshift =100pt,] (4,8) -- (4,4);
\draw[help lines,color=gray!160,step=10pt,very thick, xshift =100pt,dotted] (6,4) -- (6,6);
\draw[help lines,color=gray!160,step=10pt,very thick, xshift =100pt,dashdotted] (4,6) -- (6,6);
\draw[help lines,color=gray!160,step=10pt,very thick, xshift =100pt,dotted] (6,4) -- (6,8);
\draw[help lines,color=gray!160,step=10pt,very thick, xshift =100pt,dashdotted] (4,6) -- (8,6);
\node at (7.9,7) [right] {$\alpha=1$}; 
\node at (7.9,5) [right] {$\alpha=4$}; 
\node at (9.9,5) [right] {$\alpha=1$}; 
\node at (9.9,7) [right] {$\alpha=2$}; 
\end{tikzpicture}
\caption{The 4 subdomains and piecewise diffusion constants $\alpha$}\label{Figure_4_Subdomains}
\end{figure}
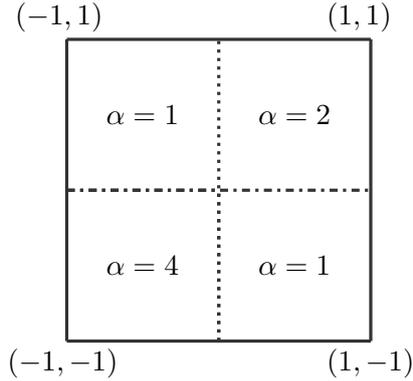

In this subsection, we test the machine learning proposed in this paper 
through the multi-material diffusion problem. 
For this aim, we consider the following interface problem:
\begin{equation}\label{interface_example_2}
\left\{
\begin{array}{rcll}
-\nabla \cdot(\alpha\nabla u)&=&f, &{\rm in}\ \bigcup\limits_{i=1}^4\Omega_i,  \\ 
u&=&0, & {\rm on}\  \partial \Omega, \\
\left[ u \right]=0,\ \ \left[\alpha\frac{\partial u}{\partial\bf n}\right]&=&g,\ \ &{\rm across}\ \Gamma,
\end{array}
\right.
\end{equation}
where the computational domain $\Omega=(-1,1)\times (-1,1)$  
consists of $4$ subdomains (see Figure \ref{Figure_4_Subdomains})
\[
\begin{aligned}
\Omega_1&: =(-1,0)\times(-1,0),\ \ \ 
\Omega_2: =(0,1)\times(-1,0),\\
\Omega_3&: =(-1,0)\times(0,1),\ \ \ \ \ 
\Omega_4: =(0,1)\times(0,1).
\end{aligned}
\]
The diffusion coefficient is a piecewise constant over the four subdomains: 
\begin{eqnarray}\label{Interface_4_coefficients}
\alpha(x_1, x_2)= \left\{
\begin{array}{ll}
\alpha_1, & (x_1, x_2) \in(-1,0] \times(-1,0], \\ 
\alpha_2, & (x_1, x_2) \in(0,1) \times(-1,0],\\ 
\alpha_3, & (x_1, x_2) \in(-1,0] \times(0,1), \\
\alpha_4, & (x_1, x_2) \in(0,1) \times(0,1).
\end{array}
\right.
\end{eqnarray}
Here, we set the exact solution is 
\begin{eqnarray}\label{Interface_Exact_2}
u(x_1, x_2)= 
\left\{
\begin{array}{ll}
c_1 \sin(k_{11}\pi x_1) \sin(k_{12}\pi x_2),  & (x_1, x_2) \in[-1,0] \times[-1,0], \\ 
c_2 \sin(k_{21}\pi x_1) \sin(k_{22}\pi x_2),  & (x_1, x_2) \in(0,1]  \times[-1,0], \\ 
c_3 \sin(k_{31}\pi x_1) \sin(k_{32}\pi x_2), & (x_1, x_2) \in[-1,0] \times(0,1], \\
c_4 \sin(k_{41}\pi x_1) \sin(k_{42}\pi x_2), & (x_1, x_2) \in(0,1]  \times(0,1],  
\end{array}
\right.
\end{eqnarray}
The source term $f(x_1, x_2)$ and the flux jump $g$
are derived from the exact solution (\ref{Interface_Exact_2}). 

First, we construct the approximation function \( u_N \). 
Consistent with the approach in the previous subsection, 
we define neural networks \( u_N^{(1)} \), \( u_N^{(2)} \), 
\( u_N^{(3)} \), and \( u_N^{(4)} \) on the four subregions. 
Since the solution of the equation belongs to \( H_0^1(\Omega) \), 
we enforce the neural networks on the subregion boundaries to be zero to 
ensure the continuity of the solution.

Next, we define neural networks \( u_N^{(5)} \), \( u_N^{(6)} \), 
\( u_N^{(7)} \), and \( u_N^{(8)} \) on the regions \( \Omega_5 := \Omega_1 \cup \Omega_2 \), 
\( \Omega_6 := \Omega_1 \cup \Omega_3 \), 
\( \Omega_7 := \Omega_2 \cup \Omega_4 \), 
and \( \Omega_8 := \Omega_3 \cup \Omega_4 \), respectively, 
and enforce them to be zero on their respective boundaries, 
ensuring that the approximation function \( u_N \) has degrees of 
freedom on the internal boundaries.

Finally, we define the neural network \( u_N^{(9)} \) 
on the region \( \Omega \) with a boundary condition of zero 
to ensure the approximation function has degrees of freedom at the 
point \( (0,0) \). 

As an example, for the neural network \( u_N^{(5)} \) 
defined on \( \Omega_5 \), we multiply it by the boundary factor function 
\( (1+x_1)(1-x_1)x_2(1+x_2) \) to obtain \( \widehat{u}_N^{(5)} \),
 ensuring that it is zero on the subregion boundaries.

Finally, we define the approximation function for the equation 
as the sum of these locally supported neural networks, i.e.,
\[
u_N = \sum_{j=1}^{9} u_N^{(j)}.
\]
It is easy to verify that \( u_N \) satisfies the Dirichlet boundary 
condition and the continuity of the solution across the interface, 
and belongs to \( H_0^1(\Omega) \).

Below, we perform numerical testing using the Ritz-type loss 
function \eqref{ritz_loss_interface}.

\subsubsection{Test for continuous interface flux}
We perform a numerical test for the case of continuous interface flux, denoted as Test (4.1). 
For this, we set the diffusion coefficients on the four subregions 
in \eqref{Interface_4_coefficients} as \( \alpha_1 = 4 \), 
\( \alpha_2 = \alpha_3 = 1 \), \( \alpha_4 = 2 \), 
and set the frequency parameters in the exact solution 
\eqref{Interface_Exact_2} as \( k_{ij} = 1 \), 
for \( i = 1, \dots, 4 \), \( j = 1, 2 \), 
and the coefficients \( c_1 = 1 \), \( c_2 = c_3 = 4 \), 
\( c_4 = 2 \). It is easy to verify that with these parameter settings, 
the flux is continuous at the interface.

We use TNNs to construct the neural network space, where each TNN's 
subnetwork is selected as an FNN with 4 hidden layers, 
each with a width of 50. The rank of the TNN is \( p = 100 \), 
and we use \( \sin(x) \) as the activation function. 
In each subregion, the TNN functions are integrated by dividing 
the 1D region into 100 subintervals, with 16 Gauss-Lobatto 
integration points per subinterval. 

For the training process, the optimization parameters are consistent 
with those used in the two-material diffusion problem. We first use 
the Adam optimizer with a learning rate of 0.001 for the first 5000 steps, 
followed by the LBFGS optimizer with a learning rate of 0.1 
for 100 steps to obtain the final result. The uniform distributed  \( 401 \times 401 \) 
points on \( [-1,1]^2 \) are chosen as the test points.

\begin{figure}[htb]
\centering
\includegraphics[width=15cm,height=4cm]{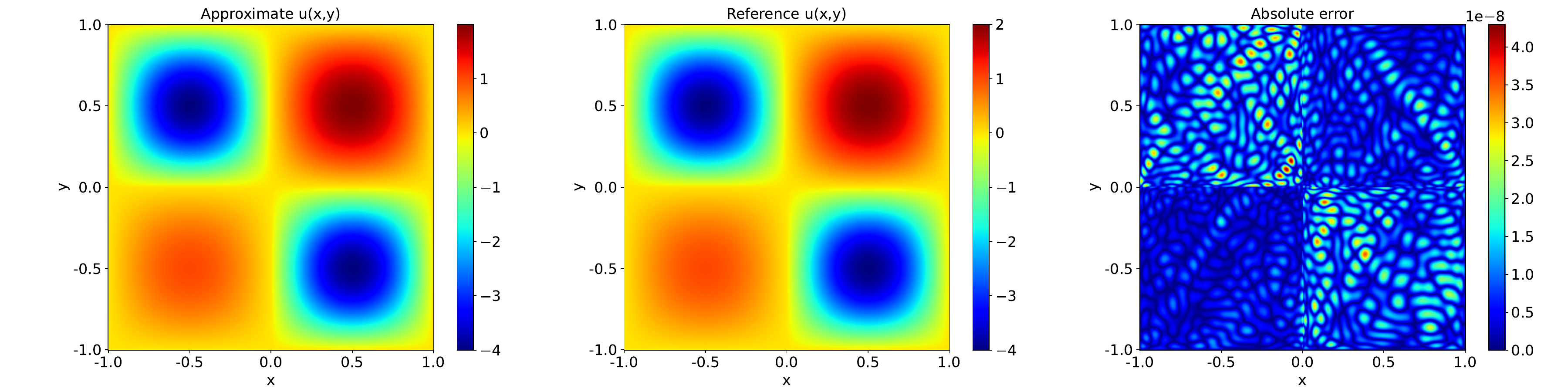}
\caption{Test (4.1): Numerical results for solving \eqref{interface_example_2} 
using the Ritz-type loss function. 
The left image is the approximate solution, the middle image 
is the exact solution, and the right image is the error distribution.}
\label{fig_interface_4_1_abs}
\end{figure}

Figure \ref{fig_interface_4_1_abs} shows the approximate solution 
and its error distribution obtained by the machine learning method, 
while Table \ref{table_interface_4_1_abs} provides the errors 
between the approximation and the exact solution. 
From these results, it can be observed that the machine 
learning method proposed in this paper can achieve  
high accuracy for this problem.
\begin{table}[htb]
\begin{center}
\caption{Numerical results of the machine learning method for Test (4.1)}
\label{table_interface_4_1_abs}
\begin{tabular}{ccc}
\hline
\( e_E \) &  \( e_{L^2} \)  &  \( e_{\text{test}} \)  \\
\hline
5.348647e-08 &   5.428252e-08 &   5.683549e-09 \\
\hline
\end{tabular}
\end{center}
\end{table}

This example is commonly used to test the performance of different machine 
learning method for solving multi-material diffusion equations. 
Table \ref{table_interface_4_1_comparison} compares the results of related 
methods based on PINNs with the method proposed in this 
paper for solving the example in this subsection. 
From the table, it can be seen that the machine learning method designed 
in this paper achieves the highest accuracy for this problem.
\begin{table}[htb]
\begin{center}
\caption{Comparison of error in solving the problem \eqref{interface_example_2} 
in Test (4.1) using different machine learning methods}\label{table_interface_4_1_comparison}
\begin{tabular}{cc}
\hline
Method &  \( e_{\text{test}} \)  \\
\hline
Standard PINN \cite{yao2023deep} &   8.96e-01 \\
\hline
nDS-PINN \cite{yao2023deep} & 5.20e-04 \\
\hline
HC-PSNN \cite{xie2024physics} & 9.00e-06 \\
\hline
Proposed Method &  5.68e-09 \\
\hline
\end{tabular}
\end{center}
\end{table}

\subsubsection{Test for jump in interface flux}
We perform a numerical test for the case of a jump in the interface flux, 
denoted as Test (4.2). 
For this, we set the diffusion coefficients on the four subregions in 
\eqref{Interface_4_coefficients} as \( \alpha_1 = 4 \), 
\( \alpha_2 = \alpha_3 = 1 \), \( \alpha_4 = 10 \), 
and set the frequency parameters in the exact solution 
\eqref{Interface_Exact_2} as \( (k_{11}, k_{12}) = (3,1) \), 
\( (k_{21}, k_{22}) = (1,2) \), \( (k_{31}, k_{32}) = (1,4) \), 
\( (k_{41}, k_{42}) = (2,2) \), 
with coefficients \( c_1 = 2 \), \( c_2 = 5 \), \( c_3 = 2 \), \( c_4 = 4 \). 
From this exact solution, it can be seen that this is an anisotropic 
problem with a flux jump at the interface.

\begin{figure}[htb]
\centering
\includegraphics[width=15cm,height=4cm]{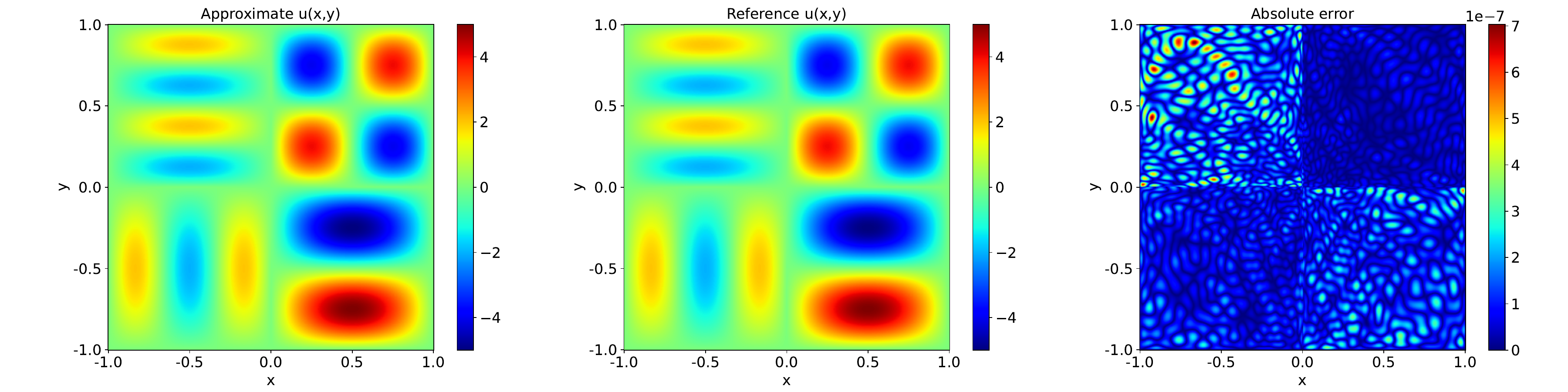}
\caption{Test (4.2): Error distribution for solving problem 
\eqref{interface_example_2} using the Ritz-type loss function.}
\label{fig_interface_4_2_abs}
\end{figure}

We use the same neural network structure, integration point configuration, 
and optimization settings as in Test (4.1). 
Figure \ref{fig_interface_4_2_abs} shows the approximate solution obtained 
by the machine learning here and the corresponding error distribution, 
while Table \ref{table_interface_4_2_abs} provides the errors 
between the approximation and the exact solution. 
From these results, it can be seen that the machine learning 
method designed in this paper also achieves high accuracy 
for this anisotropic problem with a jump in the interface flux.

\begin{figure}[htb!]
\centering
\includegraphics[width=6cm,height=6cm]{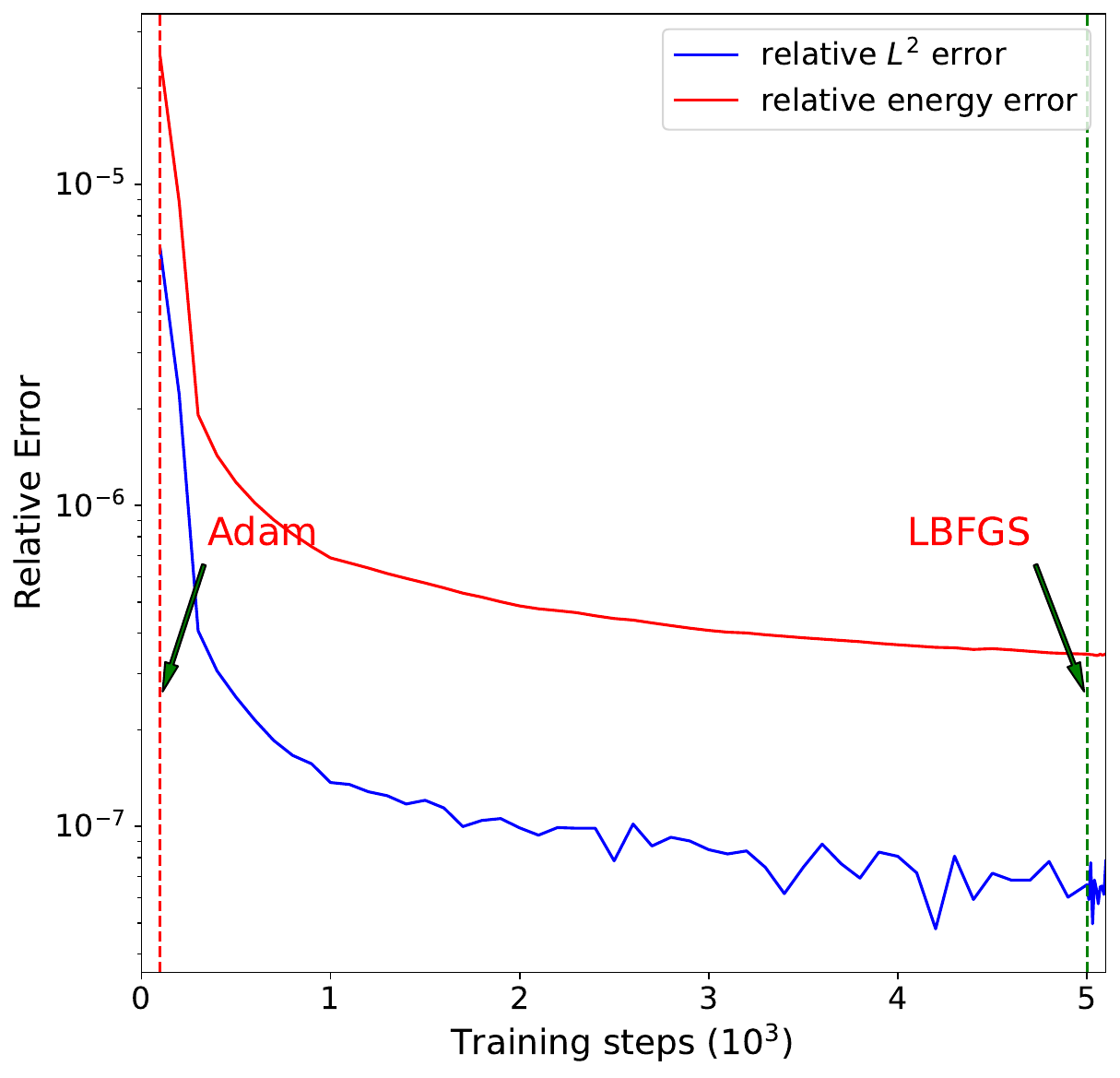}
\caption{Test (4.2): Behaviors of the relative energy error and relative \( L^2 \) 
error during the training process when solving equation \eqref{interface_example_2}.}
\label{fig_interface_4_2_iteration}
\end{figure}

Figure \ref{fig_interface_4_2_iteration} shows the behaviors  
in the relative energy error 
and the relative \( L^2 \) error during the training steps. 
From this, we can see that the training process is a relatively stable 
subspace learning process.

\begin{table}[htb]
\begin{center}
\caption{Numerical results of the machine learning method for Test (4.2)}\label{table_interface_4_2_abs}
\begin{tabular}{ccc}
\hline
\( e_E \) &  \( e_{L^2} \)  &  \( e_{\text{test}} \)  \\
\hline
3.450140e-07 &   7.798271e-08 &   6.912460e-08 \\
\hline
\end{tabular}
\end{center}
\end{table}

\subsection{A heterogeneous material problem within the domain}
In this section, we consider the case where a material is completely 
enveloped by another, which occurs commonly in heat transfer and oil reservoir simulation.
Consider the problem as follows:
\begin{eqnarray}\label{interface_circle_nonhomo}
\left\{
\begin{array}{rlll}
-\nabla \cdot(\alpha(x_1, x_2) \nabla u)&=&Q(x_1, x_2), & (x_1, x_2) 
\in \Omega_1 \cup \Omega_2, \\ 
u(x_1, \pm 2)&=&\sin \left(\frac{\pi}{4}\left(x_1^2+3\right)\right), & x_1 \in[-2,2], \\ 
u( \pm 2, x_2)&=&\sin \left(\frac{\pi}{4}\left(x_2^2+3\right)\right), & x_2 \in[-2,2], \\
\left[ u \right]=0,\ \ \left[\alpha\frac{\partial u}{\partial\bf n}\right]&=&0,\ \ &{\rm across}\ \Gamma,
\end{array}
\right.
\end{eqnarray}
where $\Omega_1:=\left\{(x_1, x_2) \mid x_1^2+x_2^2<1\right\}$, $\Omega := (-2,2)\times(-2,2)$, 
$\Omega_2:=\Omega \backslash \overline{\Omega}_1$, 
$\ \Gamma:=\left\{(x_1, x_2)\mid x_1^2+\right.$ $\left. x_2^2=1\right\}$, and
\begin{eqnarray*}
\alpha(x_1, x_2)= 
\begin{cases}
1, & (x_1, x_2) \in \Omega_1, \\ 
4, & (x_1, x_2) \in \Omega_2,
\end{cases}
\end{eqnarray*}
and the source term
\begin{eqnarray*}
Q(x_1, x_2)=\left\{
\begin{array}{c}
-4 \pi \cos \left(\pi\left(x_1^2+x_2^2-1\right)\right)
+4 \pi^2\left(x_1^2+x_2^2\right) \sin \left(\pi\left(x_1^2+x_2^2-1\right)\right), \\
(x_1, x_2) \in \Omega_1, \\
-4 \pi \cos \left(\frac{\pi}{4}\left(x_1^2+x_2^2-1\right)\right)
+\pi^2\left(x_1^2+x_2^2\right) \sin \left(\frac{\pi}{4}\left(x_1^2+x_2^2-1\right)\right), \\
(x_1, x_2) \in \Omega_2.
\end{array}
\right.
\end{eqnarray*}
Then the exact solution of this problem is
\begin{equation}
\label{interface_circle_exact_solution}
u(x_1, x_2)=
\begin{cases}
\sin \left(\pi\left(x_1^2+x_2^2-1\right)\right), & (x_1, x_2) \in \Omega_1. \\ 
\sin \left(\frac{\pi}{4}\left(x_1^2+x_2^2-1\right)\right), & (x_1, x_2) \in \Omega_2.
\end{cases}
\end{equation}

To deal with the non-homogeneous Dirichlet boundary condition, 
following the idea from finite element method, we set up 
a neural network $b_N$ to fit the boundary function $b(x_1, x_2)$ 
on $\partial \Omega$, which is defined as
\[
b(x_1, x_2) = 
\begin{cases}
\sin \left(\frac{\pi}{4}\left(x_1^2+3\right)\right), \quad 
& (x_1, x_2) \in [-2, 2] \times \{\pm 2\}, \\
\sin \left(\frac{\pi}{4}\left(x_2^2+3\right)\right), \quad 
& (x_1, x_2) \in \{\pm 2\} \times [-2, 2].
\end{cases}
\]

Then, by introducing $\varphi$ as $\varphi = u - b_N$, we transform the original problem \eqref{interface_circle_nonhomo} into the homogeneous-Dirichelet-boundary one:
\begin{eqnarray}\label{interface_circle_homo}
\left\{
\begin{array}{rlll}
-\nabla \cdot(\alpha(x_1, x_2) \nabla \varphi)&=&Q(x_1, x_2) 
+ \nabla \cdot(\alpha(x_1, x_2) \nabla b_N), & (x_1, x_2) 
\in \Omega_1 \cup \Omega_2, \\ 
\varphi(x_1, x_2) &=& 0, &(x_1, x_2) \in \partial \Omega, \\
\left[ \varphi \right]=0,\ \ \left[\alpha\frac{\partial \varphi}{\partial\bf n}\right]&=&-\left[\alpha\frac{\partial b_N}{\partial\bf n}\right]:=g,
\ \ &{\rm across}\ \Gamma.
\end{array}
\right.
\end{eqnarray}

Just as \eqref{Interface_Problem}, we define the corresponding variational 
form for \eqref{interface_circle_homo} as follows: Find $\varphi \in H_0^1(\Omega)$ such that
\[
a(\varphi, v) = (f,v), \quad \forall v \in H_0^1(\Omega),
\]
where $a(\cdot, \cdot), (\cdot, \cdot)$ are defined as follows
\[
\begin{aligned}
a(\varphi, v) &= \int_\Omega \left( \nabla v \cdot \alpha \nabla \varphi \right) d \Omega, \\ 
\quad (f, v) &= \int_\Omega Qv d\Omega 
- \int_\Omega \left( \nabla v \cdot \alpha \nabla b_N \right) d \Omega  
+ \int_\Gamma \left[\alpha\frac{\partial \varphi}{\partial\bf n}\right] v ds 
+ \int_\Gamma \left[\alpha\frac{\partial b_N}{\partial\bf n}\right] v ds \\
&= \int_\Omega Qv d\Omega - \int_\Omega \left( \nabla v \cdot \alpha \nabla b_N \right) d \Omega.
\end{aligned}
\]

As mentioned in \eqref{Interface_NN}, we define the 
trial function $\varphi_N$ by neural network functions 
$\Psi$ and $\Psi_1$, which are supported only on $\Omega$ 
and $\Omega_1$, respectively, with the following way 
\begin{equation}\label{Interface_circle_NN}
    \varphi_N = w(x_1, x_2)\Psi + \gamma(x_1, x_2)\Psi_1,
\end{equation}
where $w(x_1, x_2) = (x_1 - 2)(x_1 + 2)(x_2 - 2)(x_2 + 2)$ 
and $\gamma(x_1, x_2) = x_1^2 + x_2^2 - 1$.

Similarly, the Ritz-type loss function and the a posteriori error 
estimator-based loss function can be defined as
\begin{eqnarray}\label{interface_circle_loss}
\mathcal L_R(\varphi_N) &:=& \frac{1}{2}a(\varphi_N,\varphi_N)-(f,\varphi_N), \nonumber \\
\mathcal L_P(\varphi_N) &:=& \|\nabla\cdot(\alpha \nabla \varphi_N) 
+ Q + \nabla\cdot(\alpha \nabla b_N) \|_0^2 
+ \|g- [\mathbf n\cdot(\alpha \nabla \varphi_N)]\|_{0,\Gamma}^2,
\end{eqnarray}
respectively. We can then design a NN-based machine 
learning method utilizing Algorithm \ref{Algorithm_Subspace} 
with the loss function \eqref{interface_circle_loss} to solve \eqref{interface_circle_homo}.  
The final approximate solution for \eqref{interface_circle_nonhomo} 
is given by $u_N:= \varphi_N + b_N$. 

To interpolate the boundary function \( b(x_1,x_2) \) on \( \partial \Omega \), 
we use 100 FNNs to construct the neural network space. 
Each FNN consists of three hidden layers with a width 
of 50 and employs \( \sin(x) \) as the activation function. 
The relative \( L^2(\partial \Omega) \) error between \( b_N \) and \( b(x_1,x_2) \),
\begin{equation}\label{interface_boundary_interpolation_loss}
\frac{\|b_N(x_1,x_2) - b(x_1,x_2)\|_{L^2(\partial \Omega)}}{\|b(x_1,x_2)\|_{L^2(\partial \Omega)}},
\end{equation}
is used as the loss function. 
For the boundary integration in \eqref{interface_boundary_interpolation_loss} 
over the 4 edges of \( \partial \Omega \), the one-dimensional region is 
divided into 100 subintervals, and 16 Gauss-Lobatto quadrature points are employed. 
The Adam optimizer with a learning rate of 0.001 and a multi-step learning rate 
scheduler is utilized for 50,000 training steps. 
The final relative \( L^2(\partial \Omega) \) 
error between \( b_N \) and \( b(x_1,x_2) \) is  4.162097e-07.


To solve \eqref{interface_circle_homo} with homogeneous Dirichlet boundary conditions, 
we employ 100 FNNs to construct the neural network subspaces. 
The corresponding neural network functions $\Psi$ and $\Psi_1$ 
are used to define the trial function in \eqref{Interface_circle_NN}. 
Each FNN consists of three hidden layers with a width of 50 
and utilizes $\sin(x)$ as the activation function. 
The a posteriori error estimator-based loss function \eqref{interface_circle_loss} 
is employed to update the neural network subspace.


For the two-dimensional integration in the loss function \eqref{interface_circle_loss}, 
we design different quadrature schemes for computing the integrals over \( \Omega_1 \), 
\( \Omega_2 \), and \( \Gamma \). Let \( \mathcal{F}(\varphi_N) \) denote the integrand. 
We use a polar coordinate transformation to compute the integrals over \( \Omega_1 \) and \( \Gamma \):  
\[
\begin{aligned}
\int_{\Omega_1} \mathcal{F}(\varphi_N) d \Omega_1 &= \int_{0}^{1}\int_{0}^{2\pi} \mathcal{F}(\varphi_N) rdrd\theta, \\
\int_{\Gamma} \mathcal{F}(\varphi_N) ds &= \int_{0}^{2\pi} \mathcal{F}(\varphi_N) d\theta.
\end{aligned}
\]
For the integration with respect to \( r \), the interval \([0,1]\) is divided into 20 subintervals, 
and 8 Gauss-Lobatto points are used in each subinterval. 
For the integration with respect to \( \theta \), \( N=160 \) equispaced grid points 
in \([0, 2\pi]\) are used as quadrature points, 
with each having a quadrature weight of \( \frac{1}{N} \). 
Next, we transform the integration over \( \Omega_2 \) into the 
difference of integrals over \( \Omega \) and \( \Omega_1 \):
\[
\int_{\Omega_2}\mathcal{F}(\varphi_N) d \Omega_2 = \int_{\Omega}\mathcal{F}(\varphi_N) d \Omega - \int_{\Omega_1}\mathcal{F}(\varphi_N) d \Omega_1.
\]
To compute the integral over \( \Omega \), the region \([-2,2] \times [-2,2]\) 
is divided into 20 subintervals for both \( x_1 \) and \( x_2 \). 
Eight Gauss-Lobatto points are used in each subinterval. 
The two-dimensional quadrature points and their corresponding weights 
for integration over \( \Omega \) and \( \Omega_1 \) 
are then constructed as the Cartesian product of the one-dimensional quadrature points.

\begin{figure}[htb]
\centering
\includegraphics[width=15cm,height=4cm]{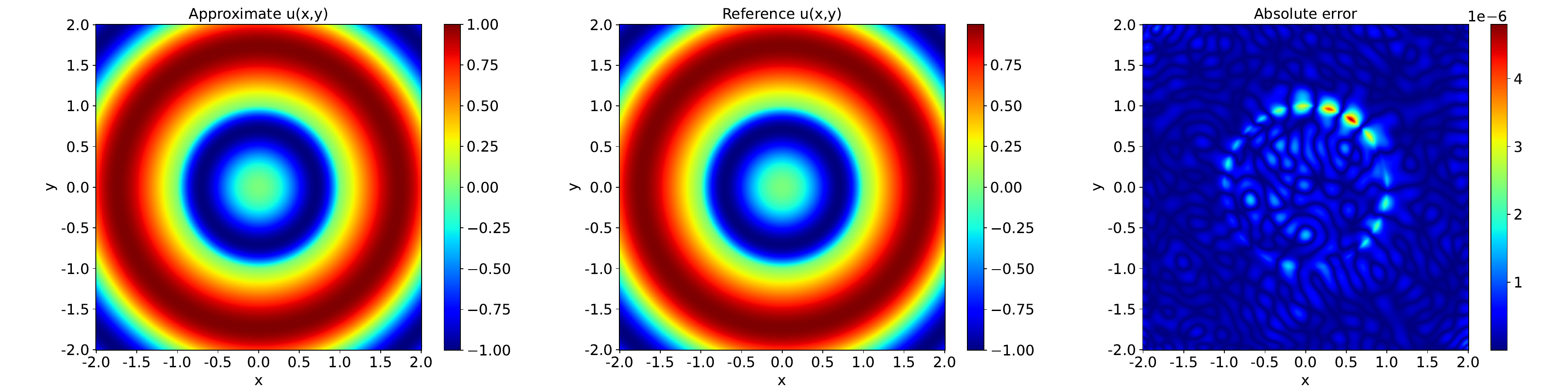}
\caption{Numerical results for solving \eqref{interface_circle_nonhomo} 
using the a posteriori error estimator-based loss function. 
The left image is the approximate solution, the middle image 
is the exact solution, and the right image is the error distribution.}
\label{fig_interface_circle}
\end{figure}

We first use the Adam optimizer with a learning rate of 0.01
for the first 5000 steps, followed by the LBFGS optimizer 
with a learning rate of 1 for 2000 steps to obtain the final result. 
The uniform distributed $401 \times 401$ points on $[-2, 2] \times [-2,2]$ 
are chosen as the test points. The relative energy norm error, relative $L^2$ norm error 
and the relative error on test points between the final approximate solution 
$u_N:=\varphi_N + b_N$ and the exact solution $u$ \eqref{interface_circle_exact_solution} 
is shown in Table \ref{table_interface_circle}. 
Figure \ref{fig_interface_circle} shows the approximate solution 
and its error distribution obtained by the proposed machine learning method.
\begin{table}[htb]
\begin{center}
\caption{Numerical results of the machine learning method for problem \eqref{interface_circle_nonhomo}}
\label{table_interface_circle}
\begin{tabular}{ccc}
\hline
\( e_E \) &  \( e_{L^2} \)  &  \( e_{\text{test}} \)  \\
\hline
3.115541e-06 & 5.595210e-07 & 5.591040e-07 \\
\hline
\end{tabular}
\end{center}
\end{table}

\section{Conclusion}
In this paper, following the idea of subspace approximation from the 
finite element method, we provide the error analysis for the neural network based 
machine learning method for solving partial differential 
equations (PDEs). The analysis reveals the significant impact of integration 
errors in the loss function computation on the accuracy of PDE solutions. 

From this understanding, we designed a type of adaptive 
neural network subspace based machine learning method 
which can achieve high accuracy for solving second order elliptic problems, 
especially those with singularities or discontinuous coefficients. 
Using the neural network subspace, we also derived the corresponding 
machine learning method and provided theoretical analysis of 
both Ritz-type and the a posterior error estimation-type loss functions, 
along with their numerical performance in solving PDEs.
Of course, the analysis and algorithms presented here can 
be easily extended to broader class of partial differential equations.

\bibliographystyle{siamplain}

\end{document}